\documentclass[12pt]{article}
\usepackage{geometry}
\usepackage[round]{natbib}
\usepackage{hyperref}
\usepackage{color}
\usepackage[tbtags]{amsmath}
\usepackage{algorithm}
\usepackage{algorithmic}
\usepackage{graphicx}
\usepackage{float}
\usepackage{subfigure}
\usepackage{cleveref}
\usepackage{authblk}
\usepackage{xcolor}
\usepackage{stmaryrd}
\usepackage{amsmath,amssymb,amsthm}

\newcommand{\X}{\mathcal{X}}

\newcommand{\y}{\mathbf{y}}

\newcommand{\w}{\mathbf{w}}
\newcommand{\B}{\mathcal{B}}

\newcommand{\bbR}{\mathbb{R}}

\newcommand{\calW}{\mathcal{W}}

\newcommand{\bB}{\mathbf{B}}
\newcommand{\bmu}{\boldsymbol{\mu}}
\newcommand{\btheta}{\boldsymbol{\theta}}
\newcommand{\bbeta}{\boldsymbol{\beta}}

\newcommand{\bepsilon}{\boldsymbol{\epsilon}}

\newcommand{\Exp}{\mathbb{E}}

\newcommand{\tp}{^\mathrm{T}}

\newcommand{\KL}{\mathrm{KL}}

\newcommand{\CV}{\mathrm{CV}}

\newcommand{\tB}{\widetilde{\B}}
\newcommand{\hB}{\widehat{\B}}

\newcommand{\tbw}{\widetilde{\B}(\mathbf{w})}

\newcommand{\sumin}{\sum_{i=1}^{n}}

\newcommand{\prodin}{\prod_{i=1}^{n}}

\newcommand{\infw}{\inf_{\mathbf{w} \in \mathcal{W}}}

\renewcommand{\vec}{\mathrm{vec}}
\renewcommand{\exp}{\mathrm{exp}}

\DeclareMathOperator{\argmin}{argmin}

\newtheorem{theorem}{Theorem} 

\newtheorem{condition}{Condition}  
\newtheorem{remark}{Remark}

\allowdisplaybreaks

\title{Improving tensor regression by optimal model averaging}
\author[1,2]{Qiushi Bu}
\author[3]{Hua Liang}
\author[1]{Xinyu Zhang}
\author[4]{Jiahui Zou}
\affil[1]{Academy of Mathematics and Systems Science, Chinese Academy of Sciences, Beijing, China}
\affil[2]{University of Chinese Academy of Sciences, Beijing, China}
\affil[3]{Department of Statistics, George Washington University, Washington, DC}
\affil[4]{School of Statistics, Capital University of Economics and Business, Beijing, China}

\begin{document}
\maketitle

\begin{abstract}
Tensors have broad applications in neuroimaging, data mining, digital marketing, etc. CANDECOMP/PARAFAC (CP) tensor decomposition can effectively reduce the number of parameters to gain dimensionality-reduction and thus plays a key role in tensor regression. However, in CP decomposition, there is uncertainty which rank to use. In this article, we develop a model averaging method to handle this uncertainty by weighting the estimators from candidate tensor regression models with different ranks. When all candidate models are misspecified, we prove that the model averaging estimator is asymptotically optimal. When correct models are included in the candidate models, we prove the consistency of parameters and the convergence of the model averaging weight. Simulations and empirical studies illustrate that the proposed method has superiority over the competition methods and has promising applications. \\

\end{abstract}

\section{Introduction} \label{sec1}
\baselineskip=20pt

In the information age, a huge amount of information can be collected, and the dimension of the data to be processed is increasing. Tensors, also called multidimensional arrays, refer to data with higher spatial dimensions, rather than high-dimensional arrays with many variables in one dimension. For example, neuroimaging data \citep{zhou2013,li2018}, longitudinal time series data \citep{wang2021,si2022}, and digital marketing data \citep{bi2018} all have 3 or more dimensions. A simple way to process a tensor is to ignore the spatial structure and vectorizing it. However, simply vectorizing a tensor and ignoring its spatial structure may result in a much larger parameter size than the sample size, and the parameters solved by this method may be inefficient \citep{yuan2016}. Taking the magnetic resonance imaging (MRI) data with dimension $\bbR^{64\times64\times64}$ as an example, the number of free parameters will be $64\times64\times64=262144$ if we use linear regression after vectorization.

Tensor decomposition solves this problem very well because it can effectively reduce the size of parameters. Various tensor decomposition methods have been developed, and two widely used methods are CP \citep{harshman1970} and Tucker tensor decomposition \citep{tucker1966}. The CP tensor decomposition can significantly compress a tensor by representing it as a linear combination of rank-1 basis tensors. An $N$-way tensor $\X\in\bbR^{I_1\times I_2\times\dots \times I_N}$ is rank-1 if it can be presented as the outer product of $N$ vectors, and the rank of a tensor is the smallest number for which a tensor can be written as the sum of rank-1 tensors \citep{kolda2009}. If we apply CP tensor decomposition with a rank-5 approximation to the same MRI data, the number of free parameters will reduce to $(64+64+64)\times 5=960$, which is much smaller than the value used in the vectorization method. \cite{guo2011} first apply the CP tensor decomposition in the linear regression models, and \cite{zhou2013} extend it to the generalized linear models (GLMs). Since then, there has been rapid development in the area of CP tensor regression. For example, \cite{lock2018} studies the case where the response variable is also a tensor, and \cite{ke2023} apply the CP tensor decomposition to quantile regressions.

Algorithms for tensor decomposition all require a given rank, but the rank is unknown in applications. \cite{haastad1989}demonstrates that calculating the rank is an NP-complete problem. The existing literature mainly uses model selection methods to choose a rank. For example, \cite{wang2022} use the Akaike information criterion (AIC) to choose rank, while \cite{zhou2013}, \cite{zhou2019} and \cite{ke2023} apply the Bayesian information criterion (BIC) to select the model. \citet{zhou2019} and \citet{ke2023} also prove the rank selected by BIC recovers the true rank for the Gaussian linear model when the true model is included in the set of candidate models. However, model selection methods are sometimes unstable since a small perturbation in the data may cause the selected model to be completely different \citep{yuan2005}. Simulations in the current paper also show the poor performance of model selection methods under high noise levels. Differing from model selection which chooses the “best” model based on certain criteria, model averaging weights the estimators of different models and avoids putting all eggs in one unevenly woven basket \citep{longford2005}.

There are two main streams of model averaging: Bayesian model averaging (BMA) and frequentist model averaging (FMA). Although BMA is flexible and can work for a wide range of models \citep{hoeting1999}, it still has drawbacks. Different prior distributions can have a profound impact on regression results, but choosing an appropriate prior distribution is often challenging and experiential. FMA can be further divided into different categories, including smoothed information criteria \citep{buckland1997,hjort2003}, adaptive weighting \citep{yuan2005,yang2001,zhang2013}, and optimal weighting \citep{wan2010,liu2015,lu2015}. Among these methods, the optimal weighting method is well developed and widely used. Initially, it focuses on linear models such as Mallows model averaging \citep{hansen2007}, jackknife model averaging \citep{hansen2012}, and heteroskedasticity robust $C_p$ \citep{liu2013}. In addition to linear regression models, researchers have studied model averaging under other frameworks such as GLMs \citep{zhang2016}, nonlinear models \citep{feng2022} and semi-parametric models \citep{li20181}. Model averaging for high-dimensional linear regression is also widely studied \citep{ando2014,yu2014,feng2020}. However, few researchers have focused on model averaging on multidimensional arrays. This article develops the optimal model averaging method for the CP tensor regression for the first time.

Based on the block relaxation algorithm proposed by \cite{zhou2013}, we develop a model averaging method based on Kullback–Leibler ($\KL$) loss to weight the estimators from candidate models with different ranks under the GLM framework. $\KL$ divergence is commonly used in evaluating the performance of the generalized linear model \citep{ando2017,zou2022} as a replacement for the square loss. To avoid overfitting, we use a $J-$fold cross-validation (CV) to calculate model averaging weights. \cite{hansen2012} and \cite{ando2017} apply leave-one-out cross-validation on linear models and generalized linear models, respectively. Once the weights are obtained, the model averaging estimator can be represented as a weighted average of  estimators from different models.

One theoretical contribution of this article is to show that our proposed estimator is asymptotically optimal when all candidate models are misspecified. Although \cite{zhou2019} and \cite{ke2023} prove that BIC model selection can choose the true rank, it only works when the candidate model set contains the true model. The asymptotic optimality guarantees that our method is at least better than BIC model selection in the case of misspecification. The simulation results also confirm this. Furthermore, the consistency of the model averaging estimator and the convergence of the weight are studied when at least one correct model is included in the set of candidate models. To the best of our knowledge, this is the first time the CP tensor regression is involved in the FMA framework. This article provides a new solution to the rank selection problem in tensor regressions and tensor decompositions.

The remainder of the article is organized as follows. Section \ref{sec2} begins with the notation and methodologies of the CP tensor regression, and then develops the model averaging procedure. Section \ref{sec3} provides theoretical results of the model averaging estimator, including asymptotic optimality under the misspecified framework, root$-n$ consistency  and the weight convergence when the set of candidate models contains the correct model. Section \ref{sec4} presents experimental results on simulated data and human brain MRI data to illustrate the advantage of the proposed method. Further discussions and conclusions are given in Section \ref{sec5}. Finally, technical proofs are provided in the Appendix.

\section{Model setup and estimation}\label{sec2}
\subsection{Notation and preliminaries}

\begin{figure}[H]
\includegraphics[width=1\linewidth]{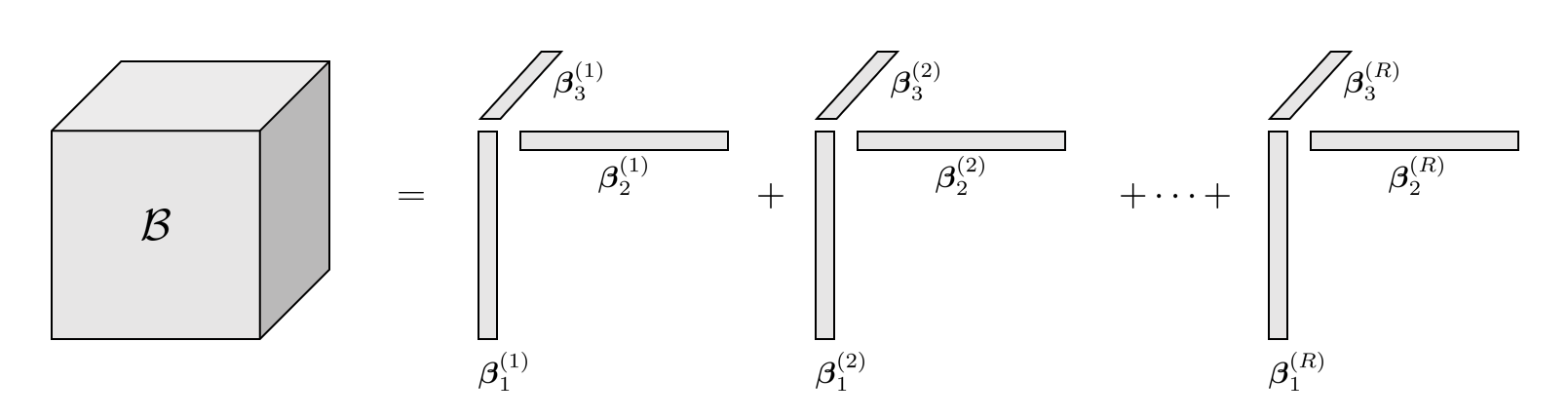}
\caption{The CP decomposition of a 3-dimensional tensor}
\label{cp3d}
\end{figure}

We introduce the notations about tensors, i.e., we use uppercase calligraphic letters, such as $\mathcal{X}$ and $\mathcal{Y}$, to represent tensors, uppercase bold letters, such as $\mathbf{X}$ and $\mathbf{Y}$, to represent matrices, and lowercase bold letters, such as $\mathbf{x}$ and $\mathbf{y}$, to represent vectors. $\left\|\cdot\right\|$ denotes the $\ell_2$ norm of a vector. The $\vec(\cdot)$ operator stacks the elements of a tensor $\B$ into a column vector. Specifically, the $(i_1, \ldots, i_D)$ entry of $\B\in \bbR^{p_1\times\dots\times p_D}$ maps to the $j$th entry of $\vec(\mathcal{B})$, where $j=1+\sum_{d=1}^D\left(i_d-1\right) \prod_{d^{\prime}=1}^{d-1} p_{d^{\prime}}$. The operator $\circ$ represents the outer product, which means each element of the outer product is the product of the corresponding vector elements: $(\bbeta^{(r)}_1\circ\cdots\circ \bbeta^{(r)}_D)_{i_1,\dots, i_D}=\beta_{1,i_1}^{(r)}\times\dots\times\beta_{D,i_D}^{(r)}$. The inner product of two tensor is written as $<\X,\B>=\vec(\X)\tp\vec(\B)$. Given a $D-$dimensional rank-$R$ tensor $\B \in \bbR^{p_1\times\dots\times p_D}$, the CP tensor decomposition factorizes it into a sum of $R$ rank-1 tensors:
$$
\B=\sum_{r=1}^{R}\bbeta_1^{(r)} \circ\cdots\circ \bbeta_D^{(r)},
$$
where $\bbeta_d^{(r)}=(\beta^{(r)}_{d,1},\dots,\beta^{(r)}_{d,p_d})\tp\in\bbR^{p_d}$ is a column vector for $d=1,\cdots,D$ and $r=1,\cdots,R$. Let $\bB_d=[\bbeta_d^{(1)},\dots,\bbeta_d^{(R)}]\in \bbR^{p_d\times R}$ for $d=1,\dots,D$, then the decomposition can be concisely expressed as $\B=\llbracket \bB_1, \ldots, \bB_D \rrbracket$ for a shorthand \citep{kolda2006}. The operator $\llbracket\cdot\rrbracket$ establishes a connection between the $D$-dimensional tensor and these $D$ matrices. Figure \ref{cp3d} gives an intuitive explanation of the CP decomposition in the 3-dimensional case.

\subsection{The CP tensor regression model}

Suppose there are $n$ independent observations $\{(y_i,\X_i), i=1,\dots,n\}$, where $\X_i\in\bbR^{p_1\times\cdots\times p_D}$ is a $D$-dimensional tensor. $p_1, \dots, p_D$ are all finite numbers and independent of $n$. We assume that the response variable $y_i$ belongs to an exponential family:
\begin{align}\label{eq:model}
\text{Pr}(y_i|\theta_i, \phi)=\exp\left\{\frac{y_i\theta_i-b(\theta_i)}{\phi}+c(y_i,\phi)\right\},
\end{align}
where $b(\cdot)$ and $c(\cdot,\cdot)$ are known functions, and $\theta_i$ and $\phi>0$ denote the natural and dispersion parameters, respectively. Under the canonical link function, $\theta_i$ relates the parameter $\B\in\bbR^{p_1\times\cdots\times p_D}$ to the covariate $\X_i$ in the form
\begin{align*}
\theta_i=<\B, \X_i>.
\end{align*}
However, the number of free parameters, $\prod_{i=1}^{D}p_i$, is often very large. Therefore, directly using the maximum likelihood method to estimate the parameters in \eqref{eq:model} will result in significant computational difficulties. If we know the rank of $\B$ is $R$, we can utilize the CP tensor decomposition to reduce the parameters, i.e.,
\begin{align*}
\B=\sum_{r=1}^{R}\bbeta_1^{(r)}\circ\cdots\circ \bbeta_D^{(r)}=\llbracket \bB_1, \ldots, \bB_D \rrbracket,
\end{align*}
where $\bbeta_d^{(r)}\in\bbR^{p_d}$ is a column vector for $d=1,\cdots,D$, $r=1,\cdots,R$, and $\bB_i=[\bbeta_i^{(1)},\dots,\bbeta_i^{(R)}]$. Under the canonical link function, we can estimate the parameters $\B$ via maximizing the loglikelihood function:
\begin{equation}
\begin{split}
&\ell\left(\bB_1, \ldots, \bB_D\right)\\
&=\sum_{i=1}^n \frac{y_i \theta_i-b\left(\theta_i\right)}{\phi}+\sum_{i=1}^n c\left(y_i, \phi\right)\\
&=\sum_{i=1}^n \frac{y_i<\B, \X_i>-b\left(<\B, \X_i>\right)}{\phi}+\sum_{i=1}^n c\left(y_i, \phi\right).
\end{split}
\label{logl}
\end{equation}
\cite{zhou2013} introduce a block relaxation algorithm to estimate the parameters, where $\mathbf{B}_d$ is sequentially updated for $d=1,\dots,D$. To address the non-uniqueness issue in the solution, they incorporate some constraints to control the scale and permutation indeterminacy for ensuring the uniqueness of estimated parameters.
\subsection{Model averaging estimation}
As mentioned above, a non-negligible problem is that we do not know which $R$ is the best since the algorithm assumes a given rank. Model averaging inspires us to average possible ranks rather than to just rely on one ``good'' $R$, which helps reduce the risk of choosing a ``poor'' $R$ by mistake. Supposing we have $S$ candidate models, the rank and the estimated parameter of the $s$th model are $r_s$ and $\widehat\B_{(s)}=\sum_{r=1}^{r_s}\widehat\bbeta_1^{(r)}\circ\cdots\circ \widehat\bbeta_D^{(r)}$, respectively. $S$ is fixed and does not change with the sample size. Furthermore, we assume that the rank of the candidate models is sorted in ascending order, i.e., $r_1<r_2<\dots<r_S$. For notational simplicity, we add $\textbf{0}$ vectors to the low-rank decomposition so that all decompositions have the same matrix size. Then
$$\widehat\B_{(s)}=\sum_{r=1}^{r_s}\widehat\bbeta_1^{(r)}\circ\cdots\circ \widehat\bbeta_D^{(r)}+\sum_{r=r_s+1}^{r_S}\boldsymbol{0}\circ\cdots\circ \boldsymbol{0}=\llbracket \widehat\bB_{(s),1}, \ldots, \widehat\bB_{(s),D} \rrbracket,
$$
where $\widehat\bB_{(s),d}=[\widehat\bbeta_d^{(1)},\dots,\widehat\bbeta_d^{(r_s)},\boldsymbol{0},\dots,\boldsymbol{0}]$ for $d=1,\dots,D$. Let $\w=(w_1\cdots,w_S)\tp\in \calW= \{\w\in[0,1]^S: \sum_{s=1}^{S}w_s=1\}$ be the weight vector, and $\widehat\bB_{d}(\w)=\sum_{s=1}^S w_s\widehat\bB_{(s),d}$ for $d=1,\dots,D$. From the perspective of model averaging, we define the model averaging estimator as the combination of several coefficient arrays with different ranks:
$$
\widehat\B(\w)=\sum_{s=1}^{S}w_s \widehat\B_{(s)}=\sum_{s=1}^{S}w_s \llbracket \widehat\bB_{(s),1}, \ldots, \widehat\bB_{(s),D} \rrbracket =\llbracket \widehat\bB_{1}(\w), \ldots, \widehat\bB_{D}(\w) \rrbracket.
$$

Let $\B_0$ be the true parameter for $\B$, $\y=(y_1,\dots,y_n)\tp$, $\bmu=\Exp(\y)$, $\theta_{0i}=<\B_0,\X_i>$ for $i=1,\dots,n$, $B_0=\sumin b(\theta_{0i})$, and $\btheta_0=(\theta_{01},\cdots,\theta_{0n})\tp$. Furthermore, let $B\{\widehat{\B}(\w)\}=\sumin b[\theta_i\{\widehat\B(\w)\}]$, and $\btheta\{\widehat{\B}(\w)\}=(\theta_1\{\widehat{\B}(\w)\}, \cdots,\theta_n\{\widehat{\B}(\w)\})\tp$. After obtaining the estimators ${\widehat\B}_{(1)},\dots,{\widehat\B}_{(S)}$, we get the $\KL$ loss of $\btheta\{\widehat{\B}(\w)\}$:
\begin{equation}
\begin{split}
\KL(\w)&=2 \Exp_{\y^*}\log\frac{\prodin \text{Pr}(y_i^*|\theta_{0i},\phi)}{\prodin \text{Pr}[y_i^*|\theta_i\{\widehat\B(\w)\},\phi]}\\
&=2\phi^{-1}B\{\widehat{\B}(\w)\}-2\phi^{-1}\bmu\tp\btheta\{\widehat{\B}(\w)\}-2\phi^{-1}B_0+2\phi^{-1}\bmu\tp\btheta_0,	
\end{split}
\label{klw}
\end{equation}
where $y_i^*$ is another realization from $\text{Pr}(\cdot |\theta_i,\phi)$ and independent of $y_i$. Note that the last two terms are not related to $\w$, so we just need to minimize $2\phi^{-1}B\{\widehat{\B}(\w)\}-2\phi^{-1}\bmu\tp\btheta\{\widehat{\B}(\w)\}$ to obtain the optimal weights.
Considering $\bmu$ is unknown, we tend to use $\y$ to replace $\bmu$, but this may lead to overfitting. Hence, we utilize the $J-$fold cross-validation to relieve this trouble and set the criterion to calculate model-averaging weights. For simplicity of expression, we assume that $n/J$ is an integer and $n_0=n/J$. We introduce additional notations for the weight selection criterion after introducing cross-validation. Let $\theta_{j,i}\{\B\}=\theta_{(j-1)n_0+i}\{\B\}$, $\X_{j,i}=\X_{(j-1)n_0+i}$, $y_{j,i}=y_{(j-1)n_0+i}$, and $\y^{(j)}=(y_{j,1},\dots,y_{j,n_0})\tp$ be the elements in the $j$th fold, and $\widehat{\B}^{[-j]}_{(s)}$ be the estimator of $\B$ of the $s$th model without the $j$th fold. Then we have $\widehat{\B}^{[-j]}(\w)=\sum_{s=1}^{S}w_s\widehat{\B}^{[-j]}_{(s)}$, $\btheta^{[-j]}\{\widehat{\B}^{[-j]}(\w)\}=(\theta_{j,1}\{\widehat{\B}^{[-j]}(\w)\}, \cdots,\theta_{j,n_0}\{\widehat{\B}^{[-j]}(\w)\})\tp$ and $B^{[-j]}\{\widehat{\B}^{[-j]}(\w)\}=\sum_{i=1}^{n_0}b[\theta_{j,i}\{\widehat{\B}^{[-j]}(\w)\}]$. \\$B\{\widehat{\B}(\w)\}$ and $\btheta\{\widehat{\B}(\w)\}$ in \eqref{klw} change into
$$
\widetilde{B}\{\tB(\w)\}=\sum_{j=1}^J B^{[-j]}\{\hB^{[-j]}(\w)\}=\sum_{j=1}^J \sum_{i=1}^{n_0} b[\theta_{j,i}\{\hB^{[-j]}(\w)\}],
$$
and
$$
\widetilde\btheta\{\tB(\w)\}=(\btheta^{[-1]}\{\widehat{\B}^{[-1]}(\w)\}\tp,\dots,\btheta^{[-J]}\{\widehat{\B}^{[-J]}(\w)\}\tp)\tp.
$$
Then the weight selection criterion $\CV_J(\w)$ is written as
\begin{equation*}
\CV_J(\w)= 2\phi^{-1}\widetilde{B}\{\tbw\}-2\phi^{-1}\y\tp\widetilde\btheta\{\tbw\}.
\end{equation*}
In general, the selection of $J$ is 5 or 10 \citep{ando2014,zhang2022, gao2023}. \cite{arlot2016} also argue that choosing $J$ to be 5 or 10 is a good choice. In this article, we assume $J=O(1)$ in the asymptotic theory part, and $J=5$ in the simulations and real data experiments.
Finally, the resultant weight vector is defined as
\begin{align}
\widehat{\w}=\argmin_{\w\in\calW}\CV_J(\w).
\label{weight}
\end{align}
Once the estimated weights are available, we can calculate the model averaging estimator, denoted by $\widehat\B(\widehat \w)=\sum_{s=1}^{S}\widehat w_s \widehat\B_{(s)}$. We term the proposed method as \textbf{T}ensor \textbf{R}egression \textbf{M}odel \textbf{A}veraging (TRMA).
\section{Asymptotic theory}\label{sec3}

\subsection{Asymptotic optimality}\label{sec3.1}

It should be emphasized that the conclusions in this paper are based on the assumption that the number of parameters $\prod_{d=1}^{D}p_d$ and the number of models $S$ are fixed, and the sample size $n$ tends to infinity. Before presenting the theorems, we first specify some regularity conditions.
Let $\bbeta_{\B}=\vec([{\bB}_{1},\dots,{\bB}_{D}])$, $\epsilon_i=y_i-\mu_i$ for $i=1,\dots,n$ and $\bepsilon=(\epsilon_1,\dots,\epsilon_n)\tp=\y-\bmu$.

\begin{condition}
For any $s=1,\dots,S$, there exists $\B^*_{(s)}=\llbracket\bB^*_{(s),1},\dots,\bB^*_{(s),D} \rrbracket $ belonging to a compact space $ \mathbb{B}$, such that
\begin{equation}
\left\|\bbeta_{\widehat\B_{(s)}}-\bbeta_{\B^*_{(s)}}\right\|=O_p(n^{-1/2}).
\label{eqcondition1}
\end{equation}

Besides, $\B^*_{(s)}$ and the true value $\B_0$ for the unknown parameter lie in the interior of the compact space $\mathbb{B}$.
\label{condition1}
\end{condition}

\begin{remark}
$\B^*_{(s)}$ is the quasi-true parameter that minimizes $\KL$ loss between the $s$th candidate model and the true model. The quasi-true parameter equals the true parameter if the working model is correctly specified.
Remark 5 in \cite{zhou2013} guarantees that the maximum likelihood estimator is the best rank-$R$ approximation of $\B_{true}$ in the sense of Kullback–Leibler distance. Given Condition \ref{condition1}, it is easy to verify that
\begin{equation}
  \left\|\vec(\widehat\B_{(s)})-\vec(\B^*_{(s)})\right\|=O_p(n^{-1/2}).
\label{eqremark1}
\end{equation}
The difference between \eqref{eqcondition1} and \eqref{eqremark1} is that \eqref{eqcondition1} requires the consistency of the vectors after the CP decomposition, while \eqref{eqremark1} indicates the consistency of the whole tensor. Equation \eqref{eqcondition1} is stronger than \eqref{eqremark1}.
\end{remark}

\begin{condition}
There is a constant $c_0$ such that
$$
\sup _{\B \in \mathbb{B}}\max _{1 \leq i \leq n}\left\|\frac{\partial\theta_i\{\B\}}{\partial \bbeta_{\B}}\right\|\leq c_0<\infty.
$$
\label{condition2}
\end{condition}

\begin{remark}
Condition \ref{condition2} states that the gradients of the systematic part are well defined, which is similar to Condition (C1) in \cite{zhang2016}.
\end{remark}
\begin{condition}
There exist positive constants $c_1$ and $c_2$ such that
$$\frac{\left\|\bmu\right\|^2}{n}\leq c_1<\infty$$
and
$$\max_{1\leq i\leq n}\left\|\vec(\X_i)\right\|\leq c_2 <\infty .$$
\label{condition3}
\end{condition}

\begin{condition}
As $n\rightarrow \infty$, there exist positive constants $c_{\min} $ and $c_{\max}$ independent of n such that $$0< c_{\min}\leq\lambda_{\min}\left\{\frac{1}{n}\sum_{i=1}^{n}\vec(\X_{i})\vec(\X_{i})\tp\right\},$$ and $$\lambda_{\max}\left\{\frac{1}{n}\sum_{i=1}^{n}\vec(\X_{i})\vec(\X_{i})\tp\right\}\leq c_{\max}<\infty,$$
where $\lambda_{\min}(\cdot)$ and $\lambda_{\max}(\cdot)$ represent the minimum and maximum eigenvalue of a matrix, respectively.
\label{condition4}
\end{condition}

\begin{remark}
Condition \ref{condition3} places restrictions on the variability of covariates, which is similar to Condition 8 in \cite{ando2014} and Conditions 2-3 in \cite{zou2022}. Condition \ref{condition4} requires that the maximum and minimum eigenvalues of the matrix $1/n\sum_{i=1}^{n}\vec(\X_{i})\vec(\X_{i})\tp$ are bounded, which is similar to Condition C.1 in \cite{zhao2020}. Condition \ref{condition4} also implies that the matrix $1/n\sum_{i=1}^{n}\vec(\X_{i})\vec(\X_{i})\tp$ is of full rank, which is a necessary condition to ensure the identifiability of the tensor regression.
\end{remark}

\begin{condition}
(i) There exists a positive constant $c_{3}$ such that for $k=1,2$,
$$
\sup _{\B \in \mathbb{B}} \max _{1 \leq i \leq n}\left|b^{(k)}\left(\theta_i\{\B\}\right)\right| \leq c_{3}<\infty ;
$$
(ii) there is a positive constant $c_{4}$ such that
$$
\inf _{\B \in \mathbb{B}}  \min _{1 \leq i \leq n} b^{''}\left(\theta_i\{\B\}\right) \geq c_{4}>0;
$$
\label{condition5}
\end{condition}

\begin{remark}
Condition \ref{condition5} regulates the behavior of function $b(\cdot )$. Common distributions such as the Poisson distribution, binary distribution and normal distribution all satisfy this condition. This condition is similar to Condition R2 in \cite{ando2017}.
\end{remark}

Let $$\KL^*(\w)=2\phi^{-1}B\{{\B}^*(\w)\}-2\phi^{-1}\bmu\tp\btheta\{{\B}^*(\w)\}-2\phi^{-1}B_0+2\phi^{-1}\bmu\tp\btheta_0,$$
and $\xi_n=\infw \KL^*(\w)$.

\begin{condition}
$n \xi_{n}^{-2}=o(1)$.
\label{condition6}
\end{condition}

\begin{remark}
Condition \ref{condition6} requires $\xi_n$ to grow faster than $n^{1/2}$, which means the candidate models are not too close to the true model. It rules out the case where the true model exists in the set of candidate models. In particular, if the ${s^*}$th model is the correct model and included in the set of candidate models, we have $\xi_n\leq \KL^*(\w_{s^*})=2\phi^{-1}(B\{{\B}^*_{(s^*)}\}-B_0)-2\phi^{-1}\bmu\tp(\btheta\{{\B}^*_{(s^*)}\}-\btheta_0)=0$, where $\w_{s^*}$ is a weight vector in which the $s^*$ element is 1 and other elements are 0. In this case, $\xi_n=0$ and $n \xi_{n}^{-2}$ will not converge to 0. This is also a common condition in model averaging articles, for example in \cite{ando2014}, \cite{liu2020} and \cite{zhang2023}.
\end{remark}

\begin{theorem}\label{thm1}
If Conditions \ref{condition1}-\ref{condition6}  hold, we have
\begin{align}
  \frac{\KL(\widehat{\w})}{\inf_{\w\in\calW}\KL(\w)}\to 1
\end{align}
in probability as $n\to\infty$.
\end{theorem}

Theorem \ref{thm1} shows that the TRMA estimator given by \eqref{weight} is asymptotically optimal in the sense of $\KL$ divergence when all candidate models are misspecified. In other words, the weight $\widehat\w$ yields a KL loss that is asymptotically equivalent to that of the infeasible optimal weight. Model misspecification is the common case because real tensor data rarely have low-rank decompositions. Theorem \ref{thm1} guarantees that our proposed method will not be worse than other methods as $n\to\infty$ under model misspecification.

\subsection{Consistency}
Following \cite{zhang2019}, the true model refers to the just-fitted model. Models $1,\dots,s_0$ are termed as the underfitted models and models $s_0+2,\dots,S$ correspond to the overfitted models. Both just-fitted models and overfitted models are called correct models.

\cite{zhou2013} prove that the estimated parameter of the true model $\bbeta_{\widehat\B}$ converges to the true parameter $\bbeta_{\B_0}$. In this article, we do not require the decomposed parameter to converge to $\bbeta_{\B_0}$, but are more concerned with the consistency of the original tensor. That is, we focus on whether the model averaging estimator converges to $\B_0$. The following theorem gives the answer.
\begin{theorem}
\label{thm2}
  If there exist correct models in the set of candidate models and Conditions \ref{condition1}-\ref{condition5} hold, then
  \begin{align}
\left\|\vec(\widehat\B(\widehat\w))-\vec(\B_0)\right\|=O_p(n^{-1/2}).
  \end{align}

\end{theorem}

Theorem \ref{thm2} states that the TRMA estimator converges to the true parameter at the rate of $\sqrt{n}$ when the correct models are contained in the candidate model set. This theorem can be viewed as a supplement of Theorem \ref{thm1}, which assumes all candidate models are misspecified.

\subsection{Weights of misclassified models}
Now we focus on the properties of weights when the correct models exist. Let
$$
\Delta_j=\left\{\left(\vec(\widehat\B^{[-j]}_{(1)})-\vec(\widehat\B^{[-j]}_{(s_0+1)})\right), \ldots,\left(\vec(\widehat\B^{[-j]}_{(s_{0})})-\vec(\widehat\B^{[-j]}_{(s_0+1)})\right)\right\}\tp,
$$
and $\widehat w_\Delta=\sum_{s=1}^{s_0}\widehat{w}_s$ be the sum of the optimal weights assigned to underfitted models.
\begin{condition}
For sufficiently large $n$, there exists a positive constant $c_{0}$ such that
$$
\min_{1\leq j\leq J}\lambda_{\min}\left({\Delta_j}{\Delta_j}\tp\right) \geq c_{0}
$$
holds almost surely.
\label{condition7}
\end{condition}
\begin{remark}
Condition \ref{condition7} is similar to Assumption 5 in \cite{liu2022}. This condition requires that there is a certain gap between the parameters of the underfitted models and those of the true model.
\end{remark}

\begin{theorem}\label{thm3}
   If Conditions \ref{condition1}-\ref{condition5} and \ref{condition7} hold and there exist correct models in the set of candidate models, we have
  $$\widehat w_\Delta=o_p(1).$$
\end{theorem}
Theorem \ref{thm3} means that the weights assigned to the underfitted models converge to zero as $n\to \infty$. In other words, the TRMA estimator ignores the ``wrong'' models and is determined only by the correct models. The three theorems take into account both misspecified and correctly specified frameworks, and in each case, our proposed method demonstrates a certain degree of superiority.

\section{Simulation Study}\label{sec4}
\subsection{2-D simulation}
\label{section2d}
\begin{figure}[H]
	\centering
	\vspace{-0.35cm}
	\subfigtopskip=2pt
	\subfigbottomskip=2pt
	\subfigcapskip=-5pt

\subfigure[Signal 1]{
		\includegraphics[width=0.2\linewidth]{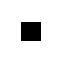}}\quad
	\subfigure[Signal 2]{
		\includegraphics[width=0.2\linewidth]{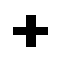}}\quad
\subfigure[Signal 3]{
		\includegraphics[width=0.2\linewidth]{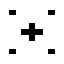}}\\
\subfigure[Signal 4]{
		\includegraphics[width=0.2\linewidth]{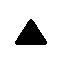}}\quad
\subfigure[Signal 5]{
		\includegraphics[width=0.2\linewidth]{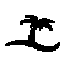}}\quad
\subfigure[Signal 6]{
		\includegraphics[width=0.2\linewidth]{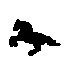}}\quad
	\caption{Six different signals, where the black area is 1 and the rest is 0}
	\label{realimg}
\end{figure}

Similar to \cite{zhou2013}, we first conduct two-dimensional shape examples. We compare the performances of different methods under a variety of signal shapes, sample sizes and noise levels. The response $y_i$ is generated from $$y_i =\eta_i+ \epsilon_i=<\X_i,\B> +\epsilon_i,$$ where $\X_i$ is a $64\times 64$ random matrix with all elements being independent standard normal distribution. The noise term $\epsilon_i \sim N(0,\sigma^2)$, and we set the noise level to be $5\%, 10\%$, and $25\% $ of the standard deviation of mean $\eta_i$ for different noise levels. $\B$ is a $64\times 64$ matrix, with the signal region being 1 and the rest 0. Figure \ref{realimg} illustrates the six signal shapes, where the black area is 1 and the rest is 0. For sample sizes, we set $n$=500, 750 and 1000. There are five candidate models with the CP tensor regression from rank-1 to rank-5, respectively. To evaluate the performance of TRMA, the following methods are used as competitors:
\begin{itemize}
\item AIC model selection \citep{akaike1973}.
\item BIC model selection \citep{schwarz1978}.
\item $\w=(w_1,\dots,w_S)\tp$, where $w_s=\exp(-\text{AIC}_s)/\sum_{i=1}^S \exp(-\text{AIC}_i)$, proposed by \cite{buckland1997}. We call it ``SAIC''.
\item $\w=(w_1,\dots,w_S)\tp$, where $w_s=\exp(-\text{BIC}_s)/\sum_{i=1}^S \exp(-\text{BIC}_i)$, which is also proposed by \cite{buckland1997}. We call it ``SBIC''.
\item $\w=(0,\dots,0,1)\tp$. This is the tensor regression model corresponding to the maximal rank decomposition. We call it ``MAX''.
\item $\w=(1/S,\dots,1/S)\tp$, where S is the number of models. This is the simple average of all models and we call it ``EQMA''.
\item Vectorize $\X$ and use LASSO to get the parameters. This model is used to illustrate the necessity of tensor decomposition on tensor data. We use 5-fold cross-validation to search for the best penalty parameter $\lambda$. We call it ``LASSO''.
\item 5-fold TRMA, the proposed method. We call it ``TRMA''.
\end{itemize}
We evaluate six different signals, where the first three signals have rank-1, rank-2, and rank-3 CP decomposition and the rest have no low-rank CP decomposition. Different experimental settings allow us to evaluate the performances of different methods when the set of candidate models is misspecified or correctly specified.

\begin{figure}[pbth]
\begin{center}
\subfigure[$n= 500$]{\label{fig1}
\includegraphics[height=.28\textheight]{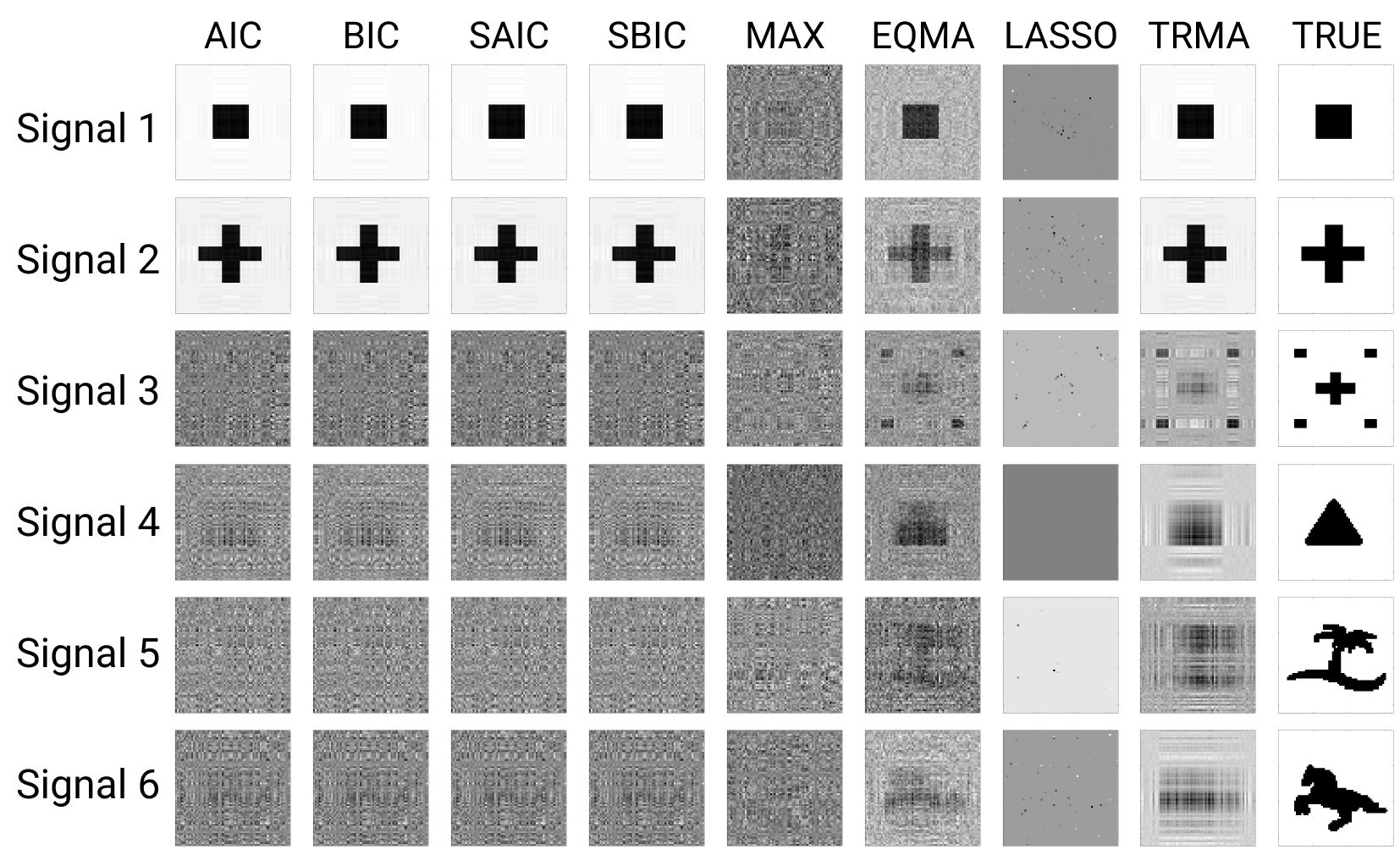}}\\
\subfigure[$n= 750$]{\label{fig2}
\includegraphics[height=.28\textheight]{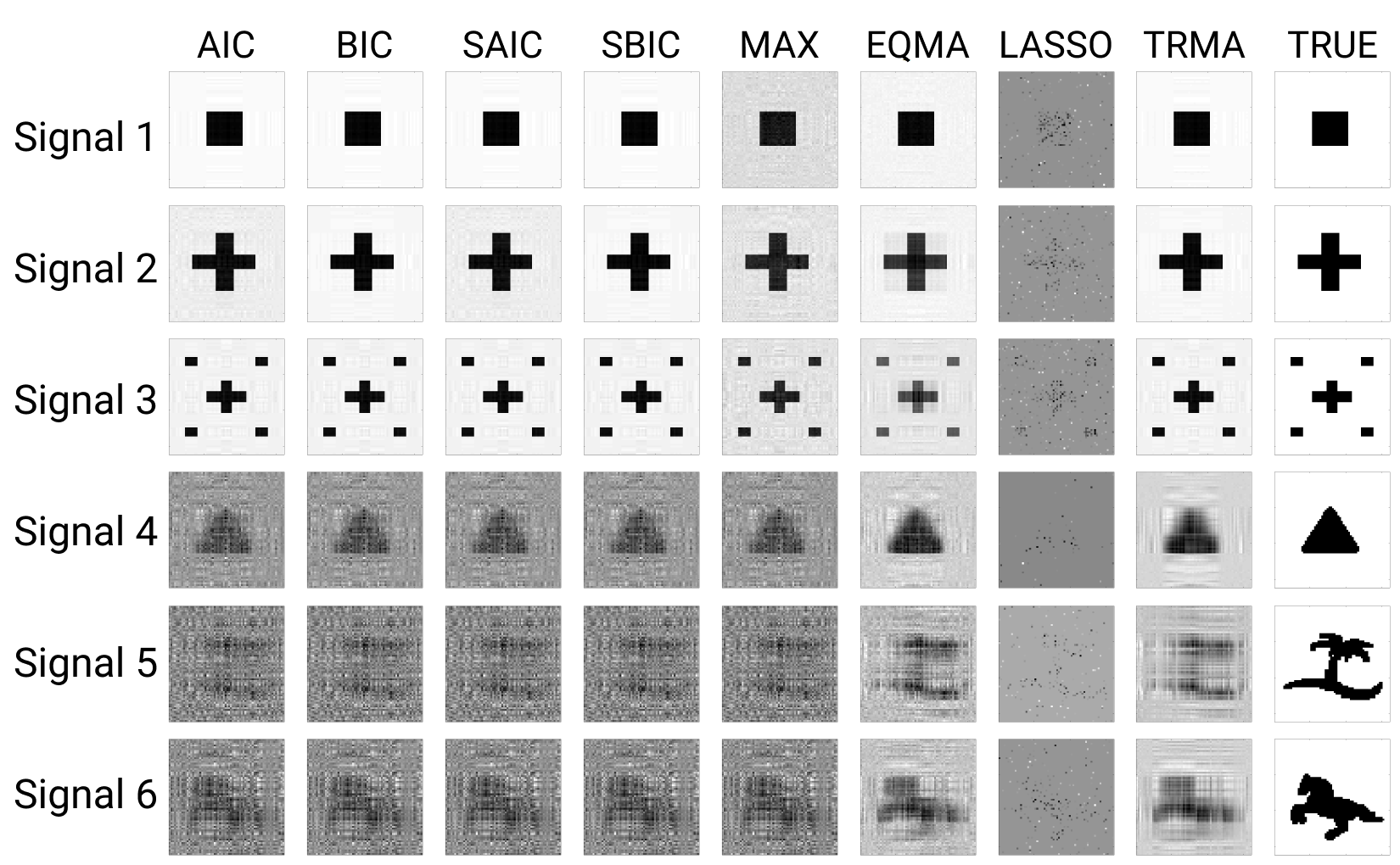}}\\
\subfigure[$n= 1000$]{\label{fig3}
\includegraphics[height=.28\textheight]{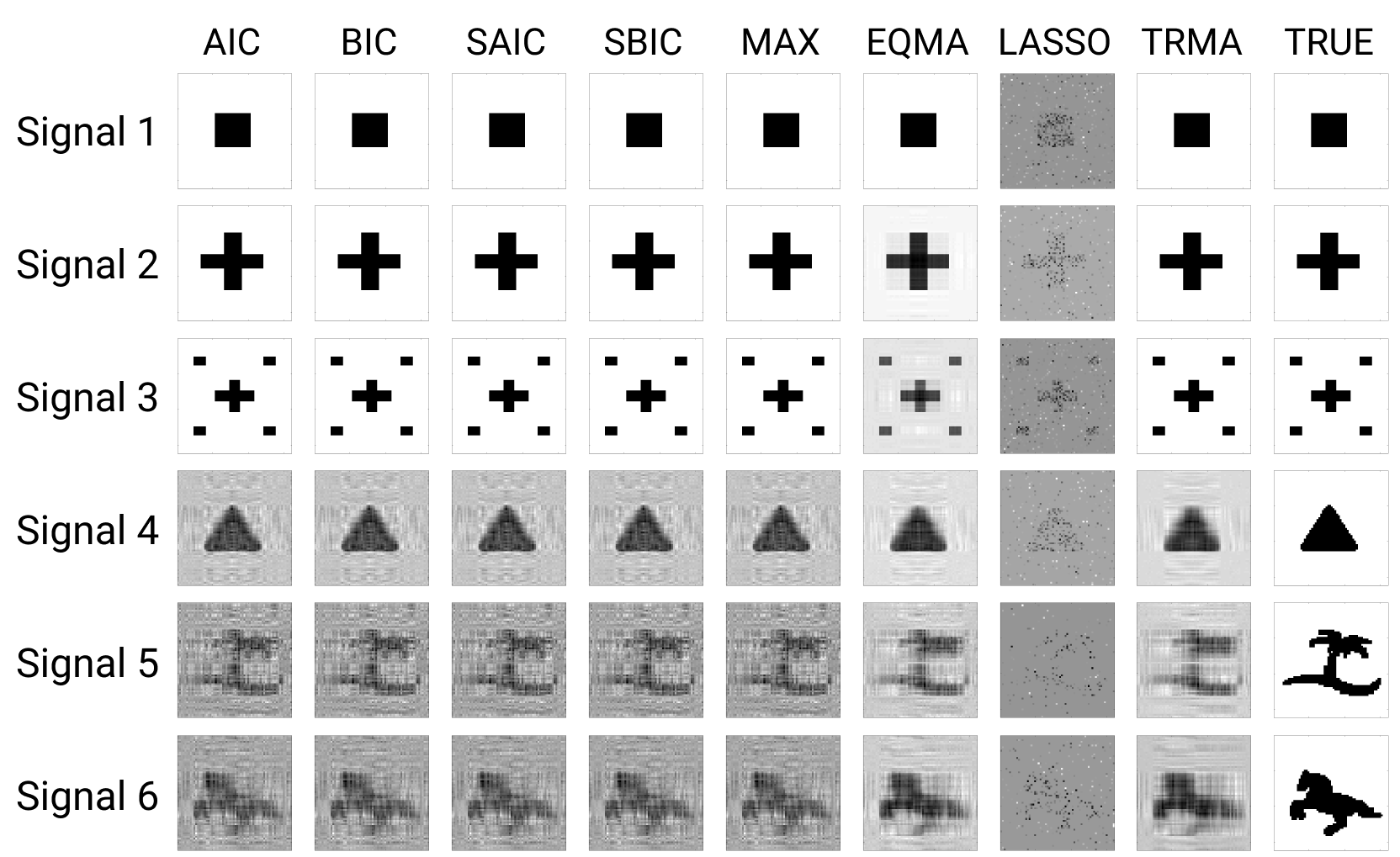}}
\caption{Recovered signal images by different methods for (a) $n=500$, (b) $750$, and (c) $1000$ and the noise level is 5\% of the standard deviation of mean $\eta_i$.}
\label{fig1-3}
\end{center}
\end{figure}

Figures \ref{fig1}-\ref{fig3} show the restoration of six signals by different methods when the noise level is 5\% of the standard deviation of mean $\eta_i$, with sample sizes of 500, 750, and 1000. For the LASSO method, the number of parameters is $64\times 64=4096$, which is much larger than the sample size. In comparison, the rank-5 CP tensor regression method has only $64\times2\times5-5^2+5=620$ parameters. The last two terms ``$-5^2+5$'' are restrictions imposed to avoid indeterminacy. As a result, in Figure \ref{fig1} we can see that the LASSO method cannot recover the signal very well due to the small sample size. Therefore, using tensor decomposition methods to reduce the number of parameters is of great significance when the number of parameters exceeds the sample size. In Figure \ref{fig1}, the recovery results of the TRMA and EQMA methods show the contour of the original image on Signals 3-6, while the results of other methods are mostly chaotic. This demonstrates the robustness of our proposed method when the sample size is relatively small. As the sample size increases from 500 to 1000, the AIC, BIC, SAIC, and SBIC methods show clearer signal recovery for Signals 3-6. However, our method still outperforms these methods by producing more distinct black and white areas with higher contrast.

The root mean squared errors (RMSEs) of $\B$ obtained from different methods are listed below to quantitatively compare the performance of different methods. Each experimental setup is replicated 100 times, and the means and standard deviations of RMSE are shown in Table \ref{tb1}. Most of the time, the TRMA method yields the smallest RMSE, with the second smallest result in the remaining cases. The LASSO method does not perform well due to the large number of parameters. When there exist correct models in the set of candidate models, i.e., in Signals 1-3, BIC and TRMA have the smallest RMSE. When the signals do not have low-rank decompositions, such as in Signals 4-6, AIC and BIC tend to choose the maximal rank model, resulting in the same RMSE as MAX. However, we can see our TRMA method outperforms even the maximal rank model. Furthermore, as the noise level increases, AIC and BIC perform worse and cannot find the true model in Signals 1-3. Under these circumstances, TRMA has a smaller RMSE, which highlights the instability of model selection under high noise levels. Although EQMA sometimes has the best performance when the signals do not have low-rank decompositions, it is not robust because it has a larger RMSE when the correct model is included in the set of candidate models. EQMA simply averages the parameters of all models. So when the true model exists, the misspecified models will greatly affect the results.

\begin{table}[pbth]
\scriptsize
\tabcolsep=0.1cm
\caption{RMSEs of $\widehat\B$ when $n=500, 750, 1000$ and the noise level is 5\%, 10\% and 20\% of the standard deviation of mean $\eta_i$.}
\label{tb1}
\begin{tabular}{cllllllllllll}
\\
\multicolumn{13}{c}{\small{$n=500$ and the noise level is 5\% of the standard deviation of mean $\eta_i$.}}\\
\hline
\multicolumn{1}{l}{} & \multicolumn{2}{c}{Signal 1} & \multicolumn{2}{c}{Signal 2} & \multicolumn{2}{c}{Signal 3} & \multicolumn{2}{c}{Signal 4} & \multicolumn{2}{c}{Signal 5} & \multicolumn{2}{c}{Signal 6} \\ \cline{2-13}
\multicolumn{1}{l}{} & Mean              & Std      & Mean              & Std      & Mean              & Std      & Mean              & Std      & Mean              & Std      & Mean              & Std      \\ \hline
AIC                  & $0.0094^{(2)}$   & 0.0022   & $0.0159^{(1)}$   & 0.0009   & 0.5596            & 0.1838   & 0.6102            & 0.0825   & 0.8796            & 0.0770   & 0.9078            & 0.0863   \\
BIC                  & $0.0091^{(1)}$   & 0.0007   & $0.0159^{(1)}$   & 0.0009   & 0.5488            & 0.1957   & 0.6102            & 0.0825   & 0.8796            & 0.0770   & 0.9103            & 0.0869   \\
SAIC                 & $0.0094^{(2)}$   & 0.0019   & $0.0159^{(1)}$   & 0.0009   & 0.5596            & 0.1838   & 0.6102            & 0.0825   & 0.8796            & 0.0770   & 0.9078            & 0.0863   \\
SBIC                 & $0.0091^{(1)}$   & 0.0007   & $0.0159^{(1)}$   & 0.0009   & 0.5488            & 0.1957   & 0.6102            & 0.0825   & 0.8796            & 0.0770   & 0.9103            & 0.0869   \\
MAX                  & 0.6083            & 0.1000   & 0.7558            & 0.0670   & 0.6392            & 0.0458   & 0.7317            & 0.0774   & 0.9057            & 0.0580   & 0.9425            & 0.0732   \\
EQMA                 & 0.1313            & 0.0218   & 0.2267            & 0.0181   & $0.2502^{(2)}$   & 0.0170   & $0.2474^{(2)}$   & 0.0139   & $0.3893^{(2)}$   & 0.0159   & $0.3876^{(2)}$   & 0.0183   \\
LASSO                & 0.3008            & 0.0030   & 0.3545            & 0.0019   & 0.2862            & 0.0030   & 0.3457            & 0.0019   & 0.4083            & 0.0014   & 0.4304            & 0.0016   \\
TRMA                 & $0.0091^{(1)}$   & 0.0007   & $0.0260^{(2)}$   & 0.0246   & $0.1975^{(1)}$   & 0.0177   & $0.1602^{(1)}$   & 0.0059   & $0.3094^{(1)}$   & 0.0079   & $0.2986^{(1)}$   & 0.0072   \\ \hline
\\
\multicolumn{13}{c}{\small{$n=750$ and the noise level is 5\% of the standard deviation of mean $\eta_i$.}}\\
\hline
\multicolumn{1}{l}{} & \multicolumn{2}{c}{Signal 1}                       & \multicolumn{2}{c}{Signal 2}                       & \multicolumn{2}{c}{Signal 3}                       & \multicolumn{2}{c}{Signal 4}                       & \multicolumn{2}{c}{Signal 5}                       & \multicolumn{2}{c}{Signal 6}                       \\ \cline{2-13}
\multicolumn{1}{l}{} & Mean              & Std      & Mean              & Std      & Mean              & Std      & Mean              & Std      & Mean              & Std      & Mean              & Std      \\ \hline
AIC                  & $0.0077^{(2)}$          & 0.0023                  & $0.0119^{(2)}$          & 0.0027                  & $0.0162^{(2)}$          & 0.0029                  & 0.2915                   & 0.0238                  & 0.5895                   & 0.0510                  & 0.5629                   & 0.0446                  \\
BIC                  & $0.0070^{(1)}$          & 0.0005                  & $0.0112^{(1)}$          & 0.0006                  & $0.0156^{(1)}$          & 0.0008                  & 0.2915                   & 0.0238                  & 0.5895                   & 0.0510                  & 0.5629                   & 0.0446                  \\
SAIC                 & $0.0077^{(2)}$          & 0.0023                  & $0.0118^{(2)}$          & 0.0025                  & $0.0162^{(2)}$          & 0.0029                  & 0.2915                   & 0.0238                  & 0.5895                   & 0.0510                  & 0.5629                   & 0.0446                  \\
SBIC                 & $0.0070^{(1)}$          & 0.0005                  & $0.0112^{(1)}$          & 0.0006                  & $0.0156^{(1)}$          & 0.0008                  & 0.2915                   & 0.0238                  & 0.5895                   & 0.0510                  & 0.5629                   & 0.0446                  \\
MAX                  & 0.0401                   & 0.0026                  & 0.0421                   & 0.0030                  & 0.0482                   & 0.0307                  & 0.2915                   & 0.0238                  & 0.5895                   & 0.0510                  & 0.5629                   & 0.0446                  \\
EQMA                 & 0.0137                   & 0.0006                  & 0.0411                   & 0.0010                  & 0.0659                   & 0.0042                  & $0.1276^{(2)}$          & 0.0048                  & $0.2621^{(2)}$          & 0.0083                  & $0.2457^{(2)}$          & 0.0074                  \\
LASSO                & 0.2889                   & 0.0062                  & 0.3463                   & 0.0050                  & 0.2718                   & 0.0062                  & 0.3356                   & 0.0047                  & 0.4023                   & 0.0042                  & 0.4250                   & 0.0043                  \\
TRMA                 & $0.0070^{(1)}$          & 0.0005                  & $0.0112^{(1)}$          & 0.0006                  & $0.0156^{(1)}$          & 0.0008                  & $0.1208^{(1)}$          & 0.0029                  & $0.2465^{(1)}$          & 0.0059                  & $0.2334^{(1)}$          & 0.0060                  \\ \hline

\\
\multicolumn{13}{c}{\small{$n=1000$ and the noise level is 5\% of the standard deviation of mean $\eta_i$.}}\\
\hline
\multicolumn{1}{l}{} & \multicolumn{2}{c}{Signal 1} & \multicolumn{2}{c}{Signal 2} & \multicolumn{2}{c}{Signal 3} & \multicolumn{2}{c}{Signal 4} & \multicolumn{2}{c}{Signal 5} & \multicolumn{2}{c}{Signal 6} \\ \cline{2-13}
\multicolumn{1}{l}{} & Mean              & Std      & Mean              & Std      & Mean              & Std      & Mean              & Std      & Mean              & Std      & Mean              & Std      \\ \hline
AIC                  & 0.0068            & 0.0021   & $0.0100^{(2)}$   & 0.0023   & $0.0126^{(2)}$   & 0.0018   & 0.1403            & 0.0091   & 0.3119            & 0.0198   & 0.2694            & 0.0199   \\
BIC                  & $0.0060^{(1)}$   & 0.0004   & $0.0091^{(1)}$   & 0.0005   & $0.0122^{(1)}$   & 0.0006   & 0.1403            & 0.0091   & 0.3119            & 0.0198   & 0.2694            & 0.0199   \\
SAIC                 & $0.0067^{(2)}$   & 0.0020   & $0.0100^{(2)}$   & 0.0022   & $0.0126^{(2)}$   & 0.0018   & 0.1403            & 0.0091   & 0.3119            & 0.0198   & 0.2694            & 0.0199   \\
SBIC                 & $0.0060^{(1)}$   & 0.0004   & $0.0091^{(1)}$   & 0.0005   & $0.0122^{(1)}$   & 0.0006   & 0.1403            & 0.0091   & 0.3119            & 0.0198   & 0.2694            & 0.0199   \\
MAX                  & 0.0272            & 0.0010   & 0.0264            & 0.0010   & 0.0254            & 0.0012   & 0.1403            & 0.0091   & 0.3119            & 0.0198   & 0.2694            & 0.0199   \\
EQMA                 & 0.0111            & 0.0006   & 0.0384            & 0.0007   & 0.0616            & 0.0012   & $0.0997^{(1)}$   & 0.0027   & $0.2112^{(1)}$   & 0.0057   & $0.1914^{(1)}$   & 0.0055   \\
LASSO                & 0.2617            & 0.0065   & 0.3275            & 0.0072   & 0.2391            & 0.0072   & 0.3170            & 0.0066   & 0.3877            & 0.0058   & 0.4130            & 0.0066   \\
TRMA                 & $0.0060^{(1)}$   & 0.0004   & $0.0091^{(1)}$   & 0.0005   & $0.0122^{(1)}$   & 0.0006   & $0.1013^{(2)}$   & 0.0024   & $0.2153^{(2)}$   & 0.0047   & $0.1980^{(2)}$   & 0.0053   \\ \hline
\\
\multicolumn{13}{c}{\small{$n=1000$ and the noise level is 10\% of the standard deviation of mean $\eta_i$.}}\\
\hline
\multicolumn{1}{l}{} & \multicolumn{2}{c}{Signal 1} & \multicolumn{2}{c}{Signal 2} & \multicolumn{2}{c}{Signal 3} & \multicolumn{2}{c}{Signal 4} & \multicolumn{2}{c}{Signal 5} & \multicolumn{2}{c}{Signal 6} \\ \cline{2-13}
\multicolumn{1}{l}{} & Mean              & Std      & Mean              & Std      & Mean              & Std      & Mean              & Std      & Mean              & Std      & Mean              & Std      \\ \hline
AIC                  & 0.0537            & 0.0022   & 0.0527            & 0.0022   & 0.0518            & 0.0025   & 0.1533            & 0.0095   & 0.3156            & 0.0184   & 0.2766            & 0.0203   \\
BIC                  & 0.0238            & 0.0061   & $0.0289^{(2)}$   & 0.0056   & $0.0350^{(2)}$   & 0.0061   & 0.1533            & 0.0095   & 0.3156            & 0.0184   & 0.2766            & 0.0203   \\
SAIC                 & 0.0537            & 0.0022   & 0.0527            & 0.0022   & 0.0518            & 0.0025   & 0.1533            & 0.0095   & 0.3156            & 0.0184   & 0.2766            & 0.0203   \\
SBIC                 & 0.0238            & 0.0060   & $0.0289^{(2)}$   & 0.0056   & $0.0349^{(2)}$   & 0.0061   & 0.1533            & 0.0095   & 0.3156            & 0.0184   & 0.2766            & 0.0203   \\
MAX                  & 0.0537            & 0.0022   & 0.0527            & 0.0022   & 0.0518            & 0.0025   & 0.1533            & 0.0095   & 0.3156            & 0.0184   & 0.2766            & 0.0203   \\
EQMA                 & $0.0219^{(2)}$   & 0.0011   & 0.0428            & 0.0013   & 0.0643            & 0.0013   & $0.1033^{(1)}$   & 0.0025   & $0.2113^{(1)}$   & 0.0062   & $0.1937^{(1)}$   & 0.0064   \\
LASSO                & 0.2636            & 0.0065   & 0.3285            & 0.0073   & 0.2417            & 0.0074   & 0.3179            & 0.0059   & 0.3892            & 0.0059   & 0.4140            & 0.0068   \\
TRMA                 & $0.0119^{(1)}$   & 0.0009   & $0.0181^{(1)}$   & 0.0009   & $0.0244^{(1)}$   & 0.0012   & $0.1045^{(2)}$   & 0.0024   & $0.2152^{(2)}$   & 0.0048   & $0.1995^{(2)}$   & 0.0051   \\ \hline
\\
\multicolumn{13}{c}{\small{$n=1000$ and the noise level is 25\% of the standard deviation of mean $\eta_i$.}}\\
\hline
\multicolumn{1}{l}{} & \multicolumn{2}{c}{Signal 1} & \multicolumn{2}{c}{Signal 2} & \multicolumn{2}{c}{Signal 3} & \multicolumn{2}{c}{Signal 4} & \multicolumn{2}{c}{Signal 5} & \multicolumn{2}{c}{Signal 6} \\ \cline{2-13}
\multicolumn{1}{l}{} & Mean              & Std      & Mean              & Std      & Mean              & Std      & Mean              & Std      & Mean              & Std      & Mean              & Std      \\ \hline
AIC                  & 0.1360            & 0.0051   & 0.1346            & 0.0047   & 0.1326            & 0.0066   & 0.2154            & 0.0115   & 0.3562            & 0.0211   & 0.3269            & 0.0232   \\
BIC                  & 0.1360            & 0.0051   & 0.1346            & 0.0047   & 0.1326            & 0.0066   & 0.2154            & 0.0115   & 0.3562            & 0.0211   & 0.3269            & 0.0232   \\
SAIC                 & 0.1360            & 0.0051   & 0.1346            & 0.0047   & 0.1326            & 0.0066   & 0.2154            & 0.0115   & 0.3562            & 0.0211   & 0.3269            & 0.0232   \\
SBIC                 & 0.1360            & 0.0051   & 0.1346            & 0.0047   & 0.1326            & 0.0066   & 0.2154            & 0.0115   & 0.3562            & 0.0211   & 0.3269            & 0.0232   \\
MAX                  & 0.1360            & 0.0051   & 0.1346            & 0.0047   & 0.1326            & 0.0066   & 0.2154            & 0.0115   & 0.3562            & 0.0211   & 0.3269            & 0.0232   \\
EQMA                 & $0.0550^{(2)}$   & 0.0021   & $0.0673^{(2)}$   & 0.0027   & $0.0839^{(2)}$   & 0.0030   & $0.1221^{(2)}$   & 0.0039   & $0.2224^{(1)}$   & 0.0064   & $0.2061^{(1)}$   & 0.0064   \\
LASSO                & 0.2719            & 0.0061   & 0.3328            & 0.0078   & 0.2521            & 0.0069   & 0.3219            & 0.0064   & 0.3912            & 0.0058   & 0.4150            & 0.0062   \\
TRMA                 & $0.0301^{(1)}$   & 0.0021   & $0.0455^{(1)}$   & 0.0026   & $0.0617^{(1)}$   & 0.0031   & $0.1201^{(1)}$   & 0.0035   & $0.2246^{(2)}$   & 0.0055   & $0.2093^{(2)}$   & 0.0056   \\ \hline
\end{tabular}

\footnotesize
Notes: The means and the standard deviations are obtained from 100 replications. Signals 1-3 correspond to the case where the correct models exist in the candidate model set, and Signals 4-6 correspond to the case where all candidate models are misspecified. The smallest and the second smallest results of each setting of experiments are flagged by (1) and (2), respectively.
\end{table}

Next, we analyze the performance of the maximal rank model. Due to the existence of the noise term, the maximal model is not the best as it will lead to overfitting. To verify this, we use the same method to generate testing data with 200 observations, and the prediction RMSE is shown in Table \ref{tb6}. The maximal rank model has the smallest RMSE on the training data, but it has a larger RMSE on the testing data. In contrast, TRMA has the smallest RMSE in most cases. Even when it is not optimal, the difference between TRMA and the optimal method is very small.

\begin{table}[pbth]
\footnotesize
\center
\tabcolsep=0.1cm
\caption{RMSEs of $\widehat \y$ on training data and testing data when $n=1000$ and the noise level is 25\% of the standard deviation of mean $\eta_i$.}
\label{tb6}
\begin{tabular}{cllllll}

\hline
\multicolumn{1}{l}{} & \multicolumn{2}{c}{Signal 1}      & \multicolumn{2}{c}{Signal 2}       & \multicolumn{2}{c}{Signal 3}       \\ \cline{2-7}
\multicolumn{1}{l}{} & Training           & Testing            & Training          & Testing             & Training           & Testing             \\ \hline
AIC                  & 0.9139          & 1.1002          & 0.8392          & 1.2097           & 0.7740          & 1.2954           \\
BIC                  & 0.9329          & $1.0739^{(1)}$ & 0.8612          & $1.1712^{(1)}$  & 0.7841          & $1.2777^{(1)}$  \\
SAIC                 & 0.9134          & 1.0987          & $0.8384^{(2)}$ & 1.2095           & $0.7739^{(2)}$ & 1.2954           \\
SBIC                 & 0.9329          & $1.0739^{(1)}$ & 0.8612          & $1.1712^{(1)}$  & 0.7841          & $1.2777^{(1)}$  \\
MAX                  & $0.3852^{(1)}$ & 2.0130          & $0.4226^{(1)}$ & 1.9758           & $0.4673^{(1)}$ & 1.9128           \\
EQMA                 & $0.5454^{(2)}$ & 1.2309          & 2.1282          & 2.6530           & 3.2246          & 4.0480           \\
LASSO                & 7.5535          & 16.8616         & 11.5289         & 21.0827          & 6.2145          & 15.4437          \\
TRMA                 & 0.9151          & $1.0747^{(2)}$ & 0.8472          & $1.1726^{(2)}$  & 0.7752          & $1.2779^{(2)}$  \\ \hline
                     & \multicolumn{2}{c}{Signal 4}      & \multicolumn{2}{c}{Signal 5}       & \multicolumn{2}{c}{Signal 6}       \\ \cline{2-7}
                     & Training           & Testing            & Training           & Testing             & Training           & Testing             \\ \hline
AIC                  & $2.1414^{(1)}$ & 8.9572          & $4.4436^{(1)}$ & 19.7966          & $4.0468^{(1)}$ & 17.2221          \\
BIC                  & $2.1414^{(1)}$ & 8.9572          & $4.4436^{(1)}$ & 19.7966          & $4.0468^{(1)}$ & 17.2221          \\
SAIC                 & $2.1414^{(1)}$ & 8.9572          & $4.4436^{(1)}$ & 19.7966          & $4.0468^{(1)}$ & 17.2221          \\
SBIC                 & $2.1414^{(1)}$ & 8.9572          & $4.4436^{(1)}$ & 19.7966          & $4.0468^{(1)}$ & 17.2221          \\
MAX                  & $2.1414^{(1)}$ & 8.9572          & $4.4436^{(1)}$ & 19.7966          & $4.0468^{(1)}$ & 17.2221          \\
EQMA                 & $3.7510^{(2)}$ & $6.4427^{(1)}$ & $7.3754^{(2)}$ & $13.4480^{(1)}$ & $7.0013^{(2)}$ & $12.3670^{(1)}$ \\
LASSO                & 10.9486         & 20.2374         & 14.9444         & 25.0418          & 17.0848         & 26.4597          \\
TRMA                 & 4.2672          & $6.5670^{(2)}$ & 9.0910          & $13.7247^{(2)}$ & 8.4024          & $12.8091^{(2)}$ \\ \hline
\end{tabular}

\footnotesize
\raggedright
Notes: The table displays $\left\|\widehat\y-\y\right\|$ on training data and testing data. The means and the standard deviations of the RMSEs are obtained from 100 replications. The smallest and the second smallest results of each setting of experiments are flagged by (1) and (2), respectively.
\end{table}

Next, we demonstrate the performance of different models under the $\KL$ loss. Table \ref{2dkl} presents the $\KL$ loss of different methods under different signals and different sample sizes. As the sample size $n$ increases, the $\KL$ loss of each method significantly decreases. When the true model exists, as in Signals 1-3, the $\KL$ loss of each method is relatively small. However, in the case of model misspecification, as in Signals 4-6, the $\KL$ loss is relatively large. After comparing different model selection and model averaging methods, it is clear that the TRMA method is quite stable, consistently producing either the best or second-best results.

\begin{table}[H]
\footnotesize
\center
\tabcolsep=0.1cm
\caption{$\KL$ loss of different methods under different signals when the noise level is 5\% of the standard deviation of mean $\eta_i$.}
\label{2dkl}
\begin{tabular}{clllllllll}
\hline
     & \multicolumn{3}{c}{Signal 1}                     & \multicolumn{3}{c}{Signal 2}                      & \multicolumn{3}{c}{Signal 3}                       \\ \cline{2-10}
     & 500             & 750            & 1000           & 500             & 750             & 1000           & 500             & 750             & 1000           \\ \hline
AIC  & 0.20            & $0.14^{(2)}$  & $0.10^{(2)}$  & $0.54^{(1)}$   & 0.31            & $0.23^{(2)}$  & 719.25          & $0.56^{(2)}$   & $0.34^{(2)}$  \\
BIC  & $0.18^{(1)}$   & $0.11^{(1)}$  & $0.07^{(1)}$  & $0.54^{(1)}$   & $0.25^{(1)}$   & $0.18^{(1)}$  & 704.99          & $0.50^{(1)}$   & $0.31^{(1)}$  \\
SAIC & $0.20^{(2)}$   & $0.14^{(2)}$  & $0.10^{(2)}$  & $0.54^{(1)}$   & $0.30^{(2)}$   & $0.23^{(2)}$  & 719.25          & $0.56^{(2)}$   & $0.34^{(2)}$  \\
SBIC & $0.18^{(1)}$   & $0.11^{(1)}$  & $0.07^{(1)}$  & $0.54^{(1)}$   & $0.25^{(1)}$   & $0.18^{(1)}$  & 704.99          & $0.50^{(1)}$   & $0.31^{(1)}$  \\
MAX  & 786.86          & 3.26           & 1.53           & 1163.78         & 3.68            & 1.45           & 834.88          & 6.35            & 1.33           \\
EQMA & 36.80           & 0.39           & 0.25           & 106.27          & 3.51            & 3.02           & $129.52^{(2)}$ & 9.20            & 7.71           \\
TRMA & $0.18^{(1)}$   & $0.11^{(1)}$  & $0.07^{(1)}$  & $2.71^{(2)}$   & $0.25^{(1)}$   & $0.18^{(1)}$  & $80.92^{(1)}$  & $0.50^{(1)}$   & $0.31^{(1)}$  \\ \hline
     & \multicolumn{3}{c}{Signal 4}                     & \multicolumn{3}{c}{Signal 5}                      & \multicolumn{3}{c}{Signal 6}                       \\ \cline{2-10}
     & 500             & 750            & 1000           & 500             & 750             & 1000           & 500             & 750             & 1000           \\ \hline
AIC  & 775.95          & 177.34         & 39.84          & 1616.23         & 698.80          & 196.70         & 1718.88         & 664.16          & 148.93         \\
BIC  & 775.95          & 177.34         & 39.84          & 1616.23         & 698.80          & 196.70         & 1723.48         & 664.16          & 148.93         \\
SAIC & 775.95          & 177.34         & 39.84          & 1616.23         & 698.80          & 196.70         & 1718.88         & 664.16          & 148.93         \\
SBIC & 775.95          & 177.34         & 39.84          & 1616.23         & 698.80          & 196.70         & 1723.48         & 664.16          & 148.93         \\
MAX  & 1111.10         & 177.34         & 39.84          & 1696.11         & 698.80          & 196.70         & 1856.59         & 664.16          & 148.93         \\
EQMA & $126.39^{(2)}$ & $33.59^{(2)}$ & $20.32^{(1)}$ & $320.38^{(2)}$ & $136.97^{(2)}$ & $90.13^{(1)}$ & $313.32^{(2)}$ & $123.99^{(2)}$ & $76.18^{(1)}$ \\
TRMA & $53.25^{(1)}$  & $29.59^{(1)}$ & $21.13^{(2)}$ & $200.74^{(1)}$ & $122.28^{(1)}$ & $93.89^{(2)}$ & $185.60^{(1)}$ & $112.13^{(1)}$ & $81.75^{(2)}$ \\ \hline
\end{tabular}

\footnotesize
\raggedright
Notes: The table displays the KL divergence of different methods under various sample sizes and signals. The means and the standard deviations of the RMSEs are obtained from 100 replications. Signals 1-3 correspond to the case where the correct models exist in the candidate model set, and Signals 4-6 correspond to the case where all candidate models are misspecified. The smallest and the second smallest results of each setting of experiments are flagged by (1) and (2), respectively.
\end{table}

\begin{figure}[H]
	\centering
	\vspace{-0.35cm}
	\subfigtopskip=2pt
	\subfigbottomskip=2pt
	\subfigcapskip=-5pt

\subfigure[Signal 4]{
		\includegraphics[width=0.31\linewidth]{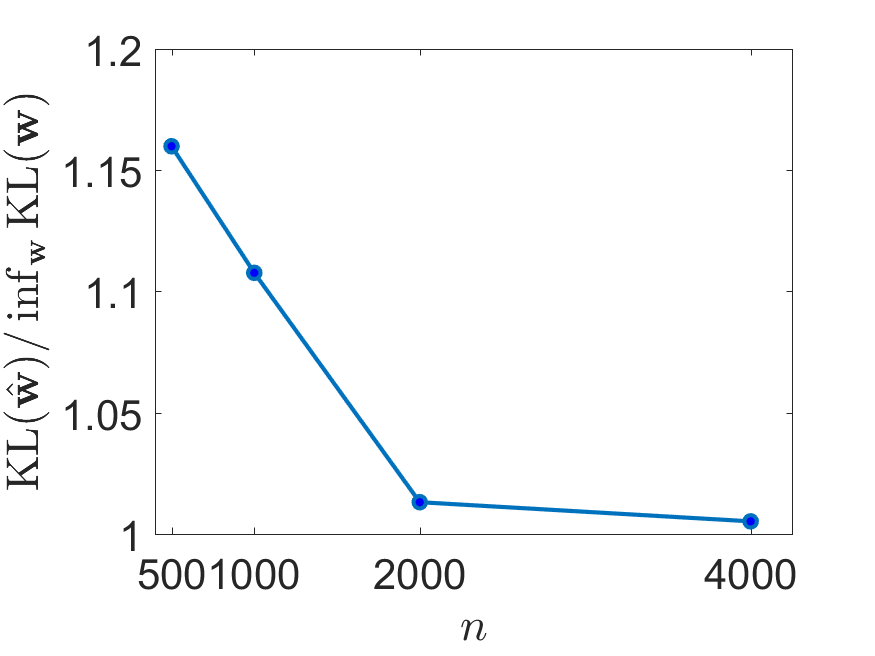}}
	\subfigure[Signal 5]{
		\includegraphics[width=0.31\linewidth]{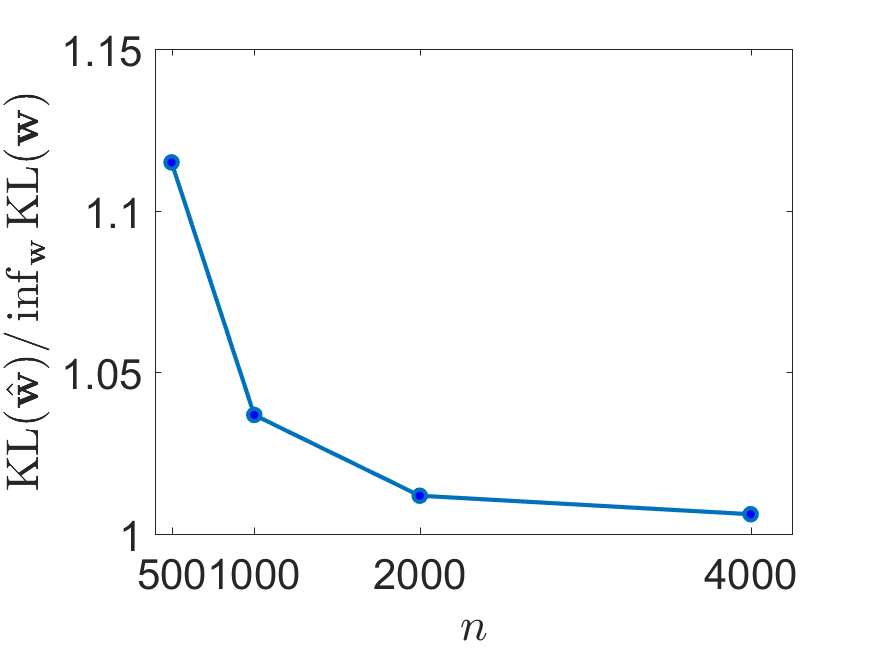}}
\subfigure[Signal 6]{
		\includegraphics[width=0.31\linewidth]{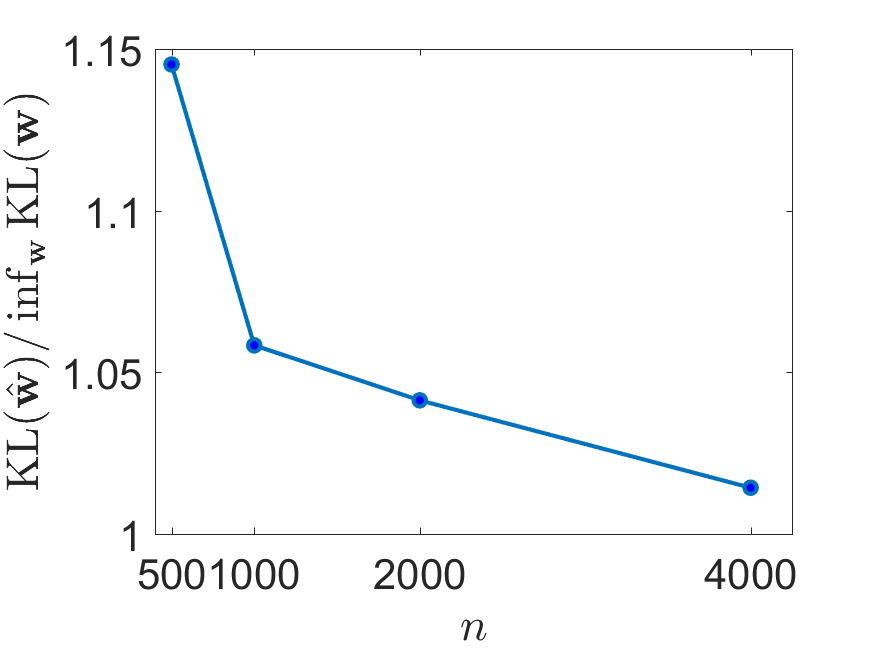}}
	\caption{The ratio of KL divergence of TRMA over the infeasible best possible model averaging.}
	\label{minkl}
\end{figure}

In Figure \ref{minkl}, the ratio of the $\KL$ loss corresponding to the model averaging estimator with respect to the smallest $\KL$ loss, i.e., $\KL(\widehat\w)/\inf_{\w}\KL(\w)$, is plotted under model misspecification. It can be seen that as $n$ increases, the ratio monotonically converges to 1, which confirms the asymptotic optimality of Theorem \ref{thm1}.

To check the root$-n$ consistency of the model averaging estimator, we set $n$ to 500, 750, 1000, 2000 and 4000 with the noise level being 5\% of the standard deviations of mean $\eta$. Table \ref{consistency} shows the values of $\left\|\vec(\widehat\B(\widehat\w))-\vec(\B_0)\right\|$ when correct models are contained. As $n$ grows, we can see that the difference between $\widehat\B(\widehat\w)$ and $\vec(\B_0)$ monotonically converges to 0, which reflects the consistency of the model averaging estimator.

\begin{table}[H]
\center
\caption{\small $\left\|\vec(\widehat\B(\widehat\w))-\vec(\B_0)\right\|$ of the TRMA estimator.}
\begin{tabular}{cccccc}
\hline
     $n$    & 500    & 750    & 1000   & 2000   & 4000   \\ \hline
Signal 1 & 0.0103 & 0.0072 & 0.0054 & 0.0041 & 0.0025 \\
Signal 2 & 0.0342 & 0.0107 & 0.0097 & 0.0059 & 0.0038 \\
Signal 3 & 0.2248 & 0.0169 & 0.0128 & 0.0073 & 0.0044 \\ \hline
\end{tabular}
\label{consistency}
\end{table}
In the experiment, we also record the weights assigned to each model to verify Theorem \ref{thm3}. It can be seen in Table \ref{wf} that the weights assigned to the underfitted models decrease as $n$ grows. When $n=500$, $\widehat w_\Delta$ is very large, most likely because of the relatively small sample size, and the parameters of each candidate model cannot be calculated well. This phenomenon gradually disappears as $n$ gets larger.
\begin{table}[H]
\center
\caption{$\widehat w_\Delta$ of the underfitted models.}
\label{wf}
\begin{tabular}{cccccc}
\hline
     $n$   & 500    & 750    & 1000   & 2000   & 4000   \\ \hline
Signal 2 & 0.0533 & 0.0023 & 0.0021 & 0.0007 & 0.0002 \\
Signal 3 & 0.8602 & 0.0112 & 0.0084 & 0.0021 & 0.0010 \\ \hline
\end{tabular}
\footnotesize

\raggedright{Notes: this table does not include Signal 1, as the CP rank of Signal 1 is 1 and no candidate model is underfitted.}
\end{table}
\subsection{3-D simulation}
Similar to Section \ref{section2d}, we evaluate normal, binomial and Poisson tensor regression on 3-D images. First we set $\eta_i=<\X_i,\B> $, where $\X_i$ and $\B\in \bbR^{32\times32\times32}$ are 3-D tensors. Then we generate the response variable through different models: for the normal model, $y_i \sim \text{Normal}(\eta_i,1)$; for the binomial model, $y_i \sim \text{Bernoulli}(p_i)$ with $p_i=\exp(0.1\eta_i)/[1+\exp(0.1\eta_i)]$; for the Poisson model, $y_i \sim \text{Poisson}(\mu_i)$ with $\mu_i=\exp(0.01\eta_i)$. We assess two 3-D signals, the first with a rank-2 CP decomposition and the second without a low-rank CP decomposition. Tables \ref{3d1}-\ref{3d2} show the means and standard deviations of the RMSE of $\widehat\B$ over 100 replications. We do not test the LASSO method because the previous experiments have shown its poor without a CP decomposition.
\begin{table}[H]
\center
\caption{RMSEs of $\widehat\B$ of 3-D Signal 1.}
\label{3d1}
\begin{tabular}{cllllll}
\hline
\multicolumn{1}{l}{} & \multicolumn{2}{c}{Normal} & \multicolumn{2}{c}{Binomial} & \multicolumn{2}{c}{Poisson} \\ \hline
\multicolumn{1}{l}{} & Mean             & Std     & Mean              & Std      & Mean              & Std     \\ \hline
AIC                  & $0.0046^{(2)}$  & 0.0007  & 1.4596            & 0.3070   & 0.4623            & 0.0894  \\
BIC                  & $0.0027^{(1)}$  & 0.0001  & $0.1174^{(1)}$   & 0.0153   & $0.3370^{(2)}$   & 0.0288  \\
SAIC                 & $0.0046^{(2)}$  & 0.0006  & 1.4595            & 0.3072   & 0.4566            & 0.0892  \\
SBIC                 & $0.0027^{(1)}$  & 0.0001  & $0.1174^{(1)}$   & 0.0153   & $0.3370^{(2)}$   & 0.0288  \\
MAX                  & 0.0072           & 0.0003  & 0.8930            & 0.1283   & 1.1687            & 0.0642  \\
EQMA                 & 0.0199           & 0.0025  & 0.4911            & 0.0402   & 0.3798            & 0.0163  \\
TRMA                 & $0.0027^{(1)}$  & 0.0002  & $0.1176^{(2)}$   & 0.0149   & $0.2874^{(1)}$   & 0.0203  \\ \hline
\end{tabular}
\footnotesize

\raggedright{Notes: The table displays $\left\|\widehat\B-\B_0\right\|$ on 3-D Signal 1. In this case, there are correct models in the set of candidate models. The means and the standard deviations of the RMSEs are obtained from 100 replications. The smallest and the second smallest results of each setting of experiments are flagged by (1) and (2), respectively.}
\end{table}

\begin{table}[H]
\center
\caption{RMSEs of $\widehat\B$ of 3-D Signal 2.}
\label{3d2}
\begin{tabular}{cllllll}
\hline
\multicolumn{1}{l}{} & \multicolumn{2}{c}{Normal} & \multicolumn{2}{c}{Binomial} & \multicolumn{2}{c}{Poisson} \\ \cline{2-7}
\multicolumn{1}{l}{} & Mean             & Std     & Mean              & Std      & Mean              & Std     \\ \hline
AIC                  & 0.1573           & 0.0275  & 1.7578            & 0.1849   & 0.5259            & 0.0932  \\
BIC                  & 0.1567           & 0.0279  & $0.1806^{(2)}$   & 0.0256   & $0.3651^{(2)}$   & 0.0251  \\
SAIC                 & 0.1573           & 0.0275  & 1.7578            & 0.1849   & 0.5143            & 0.0964  \\
SBIC                 & 0.1567           & 0.0279  & $0.1806^{(2)}$   & 0.0256   & $0.3651^{(2)}$   & 0.0251  \\
MAX                  & 0.1587           & 0.0268  & 0.9934            & 0.0910   & 1.1784            & 0.0665  \\
EQMA                 & $0.1107^{(1)}$  & 0.0120  & 0.5240            & 0.0385   & 0.4018            & 0.0182  \\
TRMA                 & $0.1170^{(2)}$  & 0.0104  & $0.1801^{(1)}$   & 0.0255   & $0.3156^{(1)}$   & 0.0187  \\ \hline
\end{tabular}
\footnotesize

\raggedright{Notes: The table displays $\left\|\widehat\B-\B_0\right\|$ on 3-D Signal 2, where all candidate models are misspecified. The means and the standard deviations of the RMSEs are obtained from 100 replications. The smallest and the second smallest results of each setting of experiments are flagged by (1) and (2), respectively.}
\end{table}
Similar to the 2-D simulation, when the signal has a low-rank CP decomposition, BIC and SBIC have the best performance, while the difference between BIC and TRMA is small. However, when the signal has no low-rank CP decomposition, TRMA performs better. From the results obtained from the binomial and Poisson models, we find that TRMA not only performs well in the normal model, but also outperforms in other distribution families. In the Poisson model, TRMA outperforms AIC and BIC when the correct models exist. Table \ref{3dkl1} also provides the $\KL$ loss of different methods. The rankings of different methods are similar to the RMSE of $\widehat\B$. Our method gives either the smallest or the second smallest $\KL$ loss.
\begin{table}[H]
\center
\caption{$\KL$ loss of 3-D Signal 1 and 3-D Signal 2.}
\label{3dkl1}
\begin{tabular}{cllllll}
\hline
\multicolumn{1}{l}{} & \multicolumn{2}{c}{Normal}          & \multicolumn{2}{c}{Binomial}      & \multicolumn{2}{c}{Poisson}       \\ \cline{2-7}
\multicolumn{1}{l}{} & Signal 1        & Signal 2          & Signal 1        & Signal 2        & Signal 1        & Signal 2        \\ \hline
AIC                  & 0.3659          & 424.4208          & 5.6251          & 9.3273          & 0.4251          & 0.6409          \\
BIC                  & $0.1192^{(1)}$ & 420.2293          & $0.1854^{(1)}$ & $0.3039^{(1)}$ & $0.2036^{(2)}$ & $0.2377^{(2)}$ \\
SAIC                 & 0.3612          & 424.4208          & 5.6341          & 9.3273          & 0.4132          & 0.5992          \\
SBIC                 & $0.1192^{(1)}$ & 420.2293          & $0.1854^{(1)}$ & $0.3039^{(1)}$ & $0.2036^{(2)}$ & $0.2377^{(2)}$ \\
MAX                  & 0.8320          & 430.3220          & 5.1145          & 5.9969          & 7.0758          & 7.8166          \\
EQMA                 & 6.3901          & $207.6599^{(1)}$ & 1.5039          & 2.3838          & 0.2576          & 0.2963          \\
TRMA                 & $0.1220^{(2)}$ & $229.8151^{(2)}$ & $0.1898^{(2)}$ & $0.3083^{(2)}$ & $0.1407^{(1)}$ & $0.1785^{(1)}$ \\ \hline
\end{tabular}
\footnotesize

\raggedright{Notes: The table displays the KL divergence of various methods under different regression model frameworks and signals. The means and the standard deviations are obtained from 100 replications. The smallest and the second smallest results of each setting of experiments are flagged by (1) and (2), respectively.}
\end{table}

Overall, the simulation results from the above experiments demonstrate that the TRMA method is effective and provides more accurate predictions compared to methods such as AIC, BIC, and SBIC. This is especially true when the data has high levels of noise and lacks low-rank decompositions.

\subsection{Real data examples}
\begin{figure}[pbth]
\begin{center}
\subfigure[50\% of the data is used as the training set]{\includegraphics[height=.23\textheight]{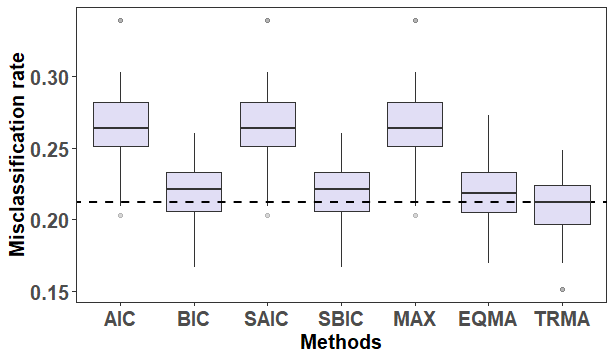}}
\subfigure[75\% of the data is used as the training set]{\includegraphics[height=.23\textheight]{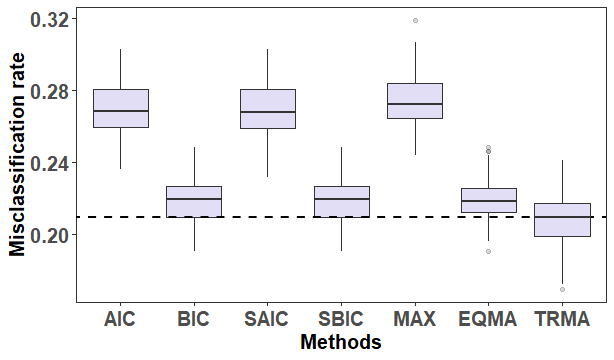}}
\subfigure[90\% of the data is used as the training set]{\includegraphics[height=.23\textheight]{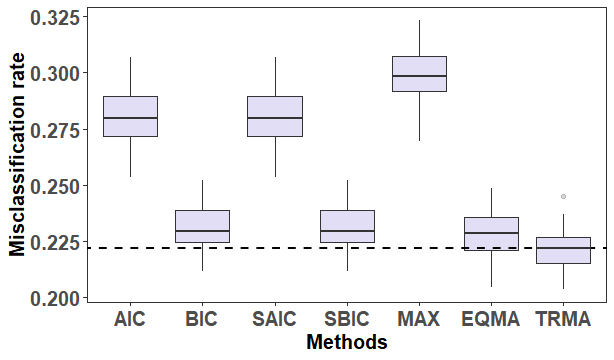}}
\caption{Misclassification rates of different methods on the skin cancer dataset when 50\%, 75\% and 90\% of the data is used as the training set. The box in the middle of the plot represents the interquartile range (IQR), which includes 50\% of the data points between the 25th percentile (Q1) and the 75th percentile (Q3). The median is represented by a horizontal line within the box. The upper and lower whiskers extending from the box indicate the maximum and minimum values within 1.5 times the IQR from the edge of the box. Any data points outside the whiskers are considered outliers and are shown as individual points beyond the whiskers. The dashed line represents the median misclassification rate of TRMA.}
\label{skin50}
\end{center}
\end{figure}

In this section, we apply the TRMA method and other methods to analyse two datasets. The first one is the skin cancer dataset (\url{https://www.isic-archive.com/}). Moles are types of skin growths or lesions, and most of them are benign and harmless. Environmental factors such as prolonged sun exposure, or changes in hormone levels such as puberty, may cause moles to darken or lesion. In this case, the benign mole can develop into melanoma, a serious skin cancer. The data contain 1800 photos of benign skin moles and 1497 photos of malignant skin moles \citep{skinref}. The purpose of this experiment is to distinguish whether the moles are benign or malignant. We downsize the origin photos to a size of $64 \times 64 \times 3$ for the convenience of calculation, and use each image as the independent variable $\X_i$. The response variable $y_i$ is binary, with 1 representing the malignant mole and 0 representing the benign mole. The set of candidate models consists of rank-1 to rank-5 CP binary tensor regression models. We then apply various model averaging and model selection methods to these five models. We conduct three sets of experiments, randomly choosing 50\%, 75\%, and 90\% of the whole data as the training set to estimate the parameters, and the remaining as the testing set to compare the misclassification rates of different models. To avoid errors caused by the number of training samples and uneven data partitioning, we repeat the experiment 100 times. Figure \ref{skin50} presents the boxplot of the misclassification rates on the testing dataset over the 100 replications.

Figure \ref{skin50} shows that the TRMA method has the highest accuracy, followed by EQMA, BIC and SBIC, and the worst performing method is the MAX method. The performance of the TRMA method is very stable and outperforms other methods under different proportions of the training set.

\begin{figure}[H]
\center
\includegraphics[width=0.7\linewidth]{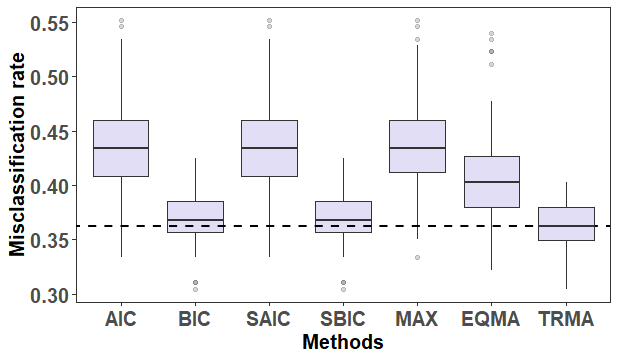}
\caption{Misclassification rates of different methods on the ADHD-200 dataset. The box in the middle of the plot represents the interquartile range (IQR), which includes 50\% of the data points between the 25th percentile (Q1) and the 75th percentile (Q3). The median is represented by a horizontal line within the box. The upper and lower whiskers extending from the box indicate the maximum and minimum values within 1.5 times the IQR from the edge of the box. Any data points outside the whiskers are considered outliers and are shown as individual points beyond the whiskers. The dashed line represents the median misclassification rate of TRMA.}
\label{adhd}
\end{figure}

The second dataset is the attention deficit hyperactivity disorder (ADHD) data (\url{http://fcon_1000.projects.nitrc.org/indi/adhd200/}). ADHD is one of the most common neurodevelopmental disorders of childhood. Common ADHD symptoms include inattention, poor impulse control, and emotional hyperactivity. The data are obtained from the ADHD-200 Global Competition datasets \citep{adhdref}, which have been split into 774 training samples and 172 testing samples. The training samples contain 285 ADHD subjects and 489 normal controls. Among them, there are 484 males and 290 females aged from 7 to 21. The raw data have been preprocessed by standard steps, including skull-stripping, segmentation with SPM12 and CAT12 toolboxes in MATLAB, etc. The size of the processed MRI images is $121 \times 145 \times 121$. We downsize each image to $10 \times 12 \times 10$ by Haar wavelet transform, and use it as covariate $\X_i$. 
The response variable $y_i$ is binary, with 1 representing ADHD and 0 representing the control group. The set of candidate models include rank-1 to rank-5 CP binary tensor regression models, and we apply various model averaging and model selection methods to these five models. Figure \ref{adhd} presents the boxplot of the misclassification rates on the testing dataset over the 100 data replications.

The TRMA method outperforms BIC and SBIC, while BIC and SBIC perform much better than the other models. The misclassification rate of the maximal rank model is the largest. This does not mean that the MRI image has a low-rank CP decomposition. The high-rank model may overfit the noise part, resulting in poor prediction results. The misclassification rate of 0.3609 is also better than that of the Tucker decomposition, the CP decomposition and the Bayes method mentioned in \cite{li2018}.

\section{Discussion}\label{sec5}
In this article, we propose a model averaging method to avoid selecting a single poor model for the CP tensor regression by weighting estimators of different models. The weights are determined by minimizing a CV-type $\KL$ divergence. We prove that the proposed estimator is asymptotically optimal when all candidate models are misspecified. When at least one candidate model is correctly specified, we also prove the TRMA estimator is root$-n$ consistent and the weights of the model averaging estimator are assigned to the correct models. The proposed method is demonstrated to be efficient and stable through both simulations and real data examples.

However, there are still some areas worthy of further research. First, while our method is based on the CP tensor regression, other tensor decomposition methods such as Tucker decomposition and tensor train decomposition also have their advantages. Combining model averaging with other tensor regression methods will be a promising avenue for future research. Second, the root$-n$ consistency of the model averaging estimator requires the number of parameters to be fixed. In practice, especially in clinical data, the sample sizes are relatively small. It is a challenging question whether the statistic has the same property after relaxing this restriction. Finally, the article only considers the case where the data are independent and identically distributed. For longitudinal data or time series data, the conclusion may be different, and this is a field that requires a further investigation.

\bibliography{ref}

\begin{thebibliography}{}

\bibitem[ADHD, 2012]{adhdref}
ADHD (2012).
\newblock {The ADHD-200 Sample}.
\newblock \url{http://fcon_1000.projects.nitrc.org/indi/adhd200/}.

\bibitem[Akaike, 1973]{akaike1973}
Akaike, H. (1973).
\newblock Maximum likelihood identification of {Gaussian} autoregressive moving
  average models.
\newblock {\em Biometrika}, 60(2):255--265.

\bibitem[Ando and Li, 2014]{ando2014}
Ando, T. and Li, K.-C. (2014).
\newblock A model-averaging approach for high-dimensional regression.
\newblock {\em Journal of the American Statistical Association},
  109(505):254--265.

\bibitem[Ando and Li, 2017]{ando2017}
Ando, T. and Li, K.-C. (2017).
\newblock A weight-relaxed model averaging approach for high-dimensional
  generalized linear models.
\newblock {\em The Annals of Statistics}, 45(6):2654--2679.

\bibitem[Arlot and Lerasle, 2016]{arlot2016}
Arlot, S. and Lerasle, M. (2016).
\newblock Choice of {$V$ for $V$-fold} cross-validation in least-squares
  density estimation.
\newblock {\em The Journal of Machine Learning Research}, 17(1):7256--7305.

\bibitem[Bi et~al., 2018]{bi2018}
Bi, X., Qu, A., and Shen, X. (2018).
\newblock Multilayer tensor factorization with applications to recommender
  systems.
\newblock {\em The Annals of Statistics}, 46(6B):3308--3333.

\bibitem[Buckland et~al., 1997]{buckland1997}
Buckland, S.~T., Burnham, K.~P., and Augustin, N.~H. (1997).
\newblock Model selection: an integral part of inference.
\newblock {\em Biometrics}, pages 603--618.

\bibitem[Feng and Liu, 2020]{feng2020}
Feng, Y. and Liu, Q. (2020).
\newblock Nested model averaging on solution path for high-dimensional linear
  regression.
\newblock {\em Stat}, 9(1):e317.

\bibitem[Feng et~al., 2022]{feng2022}
Feng, Y., Liu, Q., Yao, Q., and Zhao, G. (2022).
\newblock Model averaging for nonlinear regression models.
\newblock {\em Journal of Business \& Economic Statistics}, 40(2):785--798.

\bibitem[Gao et~al., 2023]{gao2023}
Gao, Z., Zou, J., Zhang, X., and Ma, Y. (2023).
\newblock Frequentist model averaging for envelope models.
\newblock {\em Scandinavian Journal of Statistics}.

\bibitem[Guo et~al., 2011]{guo2011}
Guo, W., Kotsia, I., and Patras, I. (2011).
\newblock Tensor learning for regression.
\newblock {\em IEEE Transactions on Image Processing}, 21(2):816--827.

\bibitem[Hansen, 2007]{hansen2007}
Hansen, B.~E. (2007).
\newblock Least squares model averaging.
\newblock {\em Econometrica}, 75(4):1175--1189.

\bibitem[Hansen and Racine, 2012]{hansen2012}
Hansen, B.~E. and Racine, J.~S. (2012).
\newblock Jackknife model averaging.
\newblock {\em Journal of Econometrics}, 167(1):38--46.

\bibitem[Harshman, 1970]{harshman1970}
Harshman, R. (1970).
\newblock Foundations of the parafac procedure: Models and conditions for an"
  explanatory" multimodal factor analysis.
\newblock {\em UCLA Working Papers in Phonetics}, 16(1):84.

\bibitem[H{\aa}stad, 1990]{haastad1989}
H{\aa}stad, J. (1990).
\newblock Tensor rank is {NP}-complete.
\newblock {\em Journal of Algorithms}, 11(4):644--654.

\bibitem[Hjort and Claeskens, 2003]{hjort2003}
Hjort, N.~L. and Claeskens, G. (2003).
\newblock Frequentist model average estimators.
\newblock {\em Journal of the American Statistical Association},
  98(464):879--899.

\bibitem[Hoeting et~al., 1999]{hoeting1999}
Hoeting, J.~A., Madigan, D., Raftery, A.~E., and Volinsky, C.~T. (1999).
\newblock Bayesian model averaging: a tutorial.
\newblock {\em Statistical Science}, 14(4):382--417.

\bibitem[ISIC, 2022]{skinref}
ISIC (2022).
\newblock Skin cancer dataset.
\newblock \url{https://www.isic-archive.com/}.

\bibitem[Ke et~al., 2023]{ke2023}
Ke, B., Zhao, W., and Wang, L. (2023).
\newblock Smoothed tensor quantile regression estimation for longitudinal data.
\newblock {\em Computational Statistics \& Data Analysis}, 178.

\bibitem[Kolda, 2006]{kolda2006}
Kolda, T. (2006).
\newblock Multilinear operators for higher-order decompositions.
\newblock {\em Sandia National Laboratories, Albuquerque, NM and Livermore,
  CA}.

\bibitem[Kolda and Bader, 2009]{kolda2009}
Kolda, T.~G. and Bader, B.~W. (2009).
\newblock Tensor decompositions and applications.
\newblock {\em SIAM Review}, 51(3):455--500.

\bibitem[Li et~al., 2018a]{li20181}
Li, J., Xia, X., Wong, W.~K., and Nott, D. (2018a).
\newblock Varying-coefficient semiparametric model averaging prediction.
\newblock {\em Biometrics}, 74(4):1417--1426.

\bibitem[Li et~al., 2018b]{li2018}
Li, X., Xu, D., Zhou, H., and Li, L. (2018b).
\newblock Tucker tensor regression and neuroimaging analysis.
\newblock {\em Statistics in Biosciences}, 10(3):520--545.

\bibitem[Liu, 2015]{liu2015}
Liu, C.-A. (2015).
\newblock Distribution theory of the least squares averaging estimator.
\newblock {\em Journal of Econometrics}, 186(1):142--159.

\bibitem[Liu and Zhang, 2022]{liu2022}
Liu, H. and Zhang, X. (2022).
\newblock Frequentist model averaging for undirected {Gaussian} graphical
  models.
\newblock {\em Biometrics}, pages 1--13.

\bibitem[Liu and Okui, 2013]{liu2013}
Liu, Q. and Okui, R. (2013).
\newblock Heteroscedasticity-robust {$C_p$} model averaging.
\newblock {\em The Econometrics Journal}, 16(3):463--472.

\bibitem[Liu et~al., 2020]{liu2020}
Liu, Q., Yao, Q., and Zhao, G. (2020).
\newblock Model averaging estimation for conditional volatility models with an
  application to stock market volatility forecast.
\newblock {\em Journal of Forecasting}, 39(5):841--863.

\bibitem[Lock, 2018]{lock2018}
Lock, E.~F. (2018).
\newblock Tensor-on-tensor regression.
\newblock {\em Journal of Computational and Graphical Statistics},
  27(3):638--647.

\bibitem[Longford, 2005]{longford2005}
Longford, N.~T. (2005).
\newblock {Model Selection and Efficiency: Is `Which Model...?' the Right
  Question?}
\newblock {\em Journal of the Royal Statistical Society. Series A (Statistics
  in Society)}, 168(3):469--472.

\bibitem[Lu and Su, 2015]{lu2015}
Lu, X. and Su, L. (2015).
\newblock Jackknife model averaging for quantile regressions.
\newblock {\em Journal of Econometrics}, 188(1):40--58.

\bibitem[Schwarz, 1978]{schwarz1978}
Schwarz, G. (1978).
\newblock Estimating the dimension of a model.
\newblock {\em The Annals of Statistics}, 6(2):461--464.

\bibitem[Si et~al., 2022]{si2022}
Si, Y., Zhang, Y., and Li, G. (2022).
\newblock An efficient tensor regression for high-dimensional data.
\newblock {\em arXiv preprint arXiv:2205.13734}.

\bibitem[Tucker, 1966]{tucker1966}
Tucker, L.~R. (1966).
\newblock Some mathematical notes on three-mode factor analysis.
\newblock {\em Psychometrika}, 31(3):279--311.

\bibitem[Wan et~al., 2010]{wan2010}
Wan, A.~T., Zhang, X., and Zou, G. (2010).
\newblock Least squares model averaging by {Mallows} criterion.
\newblock {\em Journal of Econometrics}, 156(2):277--283.

\bibitem[Wang et~al., 2021]{wang2021}
Wang, D., Zheng, Y., and Li, G. (2021).
\newblock High-dimensional low-rank tensor autoregressive time series modeling.
\newblock {\em arXiv preprint arXiv:2101.04276}.

\bibitem[Wang et~al., 2022]{wang2022}
Wang, J., Hou, J., and Eldar, Y.~C. (2022).
\newblock {Tensor robust principal component analysis from multilevel quantized
  observations}.
\newblock {\em IEEE Transactions on Information Theory}, 69(1):383--406.

\bibitem[Yang, 2001]{yang2001}
Yang, Y. (2001).
\newblock Adaptive regression by mixing.
\newblock {\em Journal of the American Statistical Association},
  96(454):574--588.

\bibitem[Yu and Feng, 2014]{yu2014}
Yu, Y. and Feng, Y. (2014).
\newblock Modified cross-validation for penalized high-dimensional linear
  regression models.
\newblock {\em Journal of Computational and Graphical Statistics},
  23(4):1009--1027.

\bibitem[Yuan and Zhang, 2016]{yuan2016}
Yuan, M. and Zhang, C.-H. (2016).
\newblock On tensor completion via nuclear norm minimization.
\newblock {\em Foundations of Computational Mathematics}, 16(4):1031--1068.

\bibitem[Yuan and Yang, 2005]{yuan2005}
Yuan, Z. and Yang, Y. (2005).
\newblock Combining linear regression models: When and how?
\newblock {\em Journal of the American Statistical Association},
  100(472):1202--1214.

\bibitem[Zhang et~al., 2019]{zhou2019}
Zhang, X., Li, L., Zhou, H., Zhou, Y., Shen, D., and ADNI (2019).
\newblock Tensor generalized estimating equations for longitudinal imaging
  analysis.
\newblock {\em Statistica Sinica}, 29(4):1977--2005.

\bibitem[Zhang and Liu, 2019]{zhang2019}
Zhang, X. and Liu, C.-A. (2019).
\newblock Inference after model averaging in linear regression models.
\newblock {\em Econometric Theory}, 35(4):816--841.

\bibitem[Zhang and Liu, 2023]{zhang2022}
Zhang, X. and Liu, C.-A. (2023).
\newblock Model averaging prediction by {$K$}-fold cross-validation.
\newblock {\em Journal of Econometrics}, 235(1):280--301.

\bibitem[Zhang et~al., 2013]{zhang2013}
Zhang, X., Lu, Z., and Zou, G. (2013).
\newblock Adaptively combined forecasting for discrete response time series.
\newblock {\em Journal of Econometrics}, 176(1):80--91.

\bibitem[Zhang et~al., 2016]{zhang2016}
Zhang, X., Yu, D., Zou, G., and Liang, H. (2016).
\newblock Optimal model averaging estimation for generalized linear models and
  generalized linear mixed-effects models.
\newblock {\em Journal of the American Statistical Association},
  111(516):1775--1790.

\bibitem[Zhang et~al., 2023]{zhang2023}
Zhang, X., Zhang, X., and Ma, Y. (2023).
\newblock A model-averaging treatment of multiple instruments in {Poisson}
  models with errors.
\newblock {\em Canadian Journal of Statistics}, 51(1):173--198.

\bibitem[Zhao et~al., 2020]{zhao2020}
Zhao, S., Liao, J., and Yu, D. (2020).
\newblock Model averaging estimator in ridge regression and its large sample
  properties.
\newblock {\em Statistical Papers}, 61(4):1719--1739.

\bibitem[Zhou et~al., 2013]{zhou2013}
Zhou, H., Li, L., and Zhu, H. (2013).
\newblock Tensor regression with applications in neuroimaging data analysis.
\newblock {\em Journal of the American Statistical Association},
  108(502):540--552.

\bibitem[Zou et~al., 2022]{zou2022}
Zou, J., Wang, W., Zhang, X., and Zou, G. (2022).
\newblock Optimal model averaging for divergent-dimensional {Poisson}
  regressions.
\newblock {\em Econometric Reviews}, 41(7):775--805.

\end{thebibliography}

\end{document}